\documentclass[final,5p,times,twocolumn]{elsarticle}

\usepackage[export]{adjustbox}
\usepackage{subfig}
\usepackage{graphics}
\usepackage{lipsum}
\usepackage{caption}
\usepackage{amsmath}
\usepackage{bm}
\usepackage{amssymb}
\usepackage[nodots,nocompress]{numcompress}
\usepackage{lineno}
\usepackage{hyperref}
\usepackage{color}
\usepackage{algorithm,algorithmic}

\usepackage[normalem]{ulem}
\newcommand{\revise}[2]{#2}

\graphicspath{{images/}}
\journal{Solid and Physical Modeling 2018}

\begin{document}

\begin{frontmatter}

\title{Continuous Optimization of Adaptive Quadtree Structures}



\author{Jun Wu}\ead{j.wu-1@tudelft.nl}
\address{Depart. of Design Engineering, Delft University of Technology, The Netherlands}

\begin{abstract}
We present a novel \textit{continuous} optimization method to the \textit{discrete} problem of quadtree optimization. The optimization aims at achieving a quadtree structure with the highest mechanical stiffness, where the edges in the quadtree are interpreted as structural elements carrying mechanical loads. We formulate quadtree optimization as a continuous material distribution problem. The discrete design variables (i.e., to refine or not to refine) are replaced by continuous variables on multiple levels in the quadtree hierarchy. In discrete quadtree optimization, a cell is only eligible for refinement if its parent cell has been refined. We propose a continuous analogue to this dependency for continuous multi-level design variables, and integrate it in the iterative optimization process. Our results show that the continuously optimized quadtree structures perform much stiffer than uniform patterns and the heuristically optimized counterparts. We demonstrate the use of adaptive structures as lightweight infill for 3D printed parts, where uniform geometric patterns have been typically used in practice.




\end{abstract}


\end{frontmatter}


\section{Introduction}\label{sec:Introduction}

\begin{figure}[t]
\centering
\small
\def\svgwidth{\linewidth}
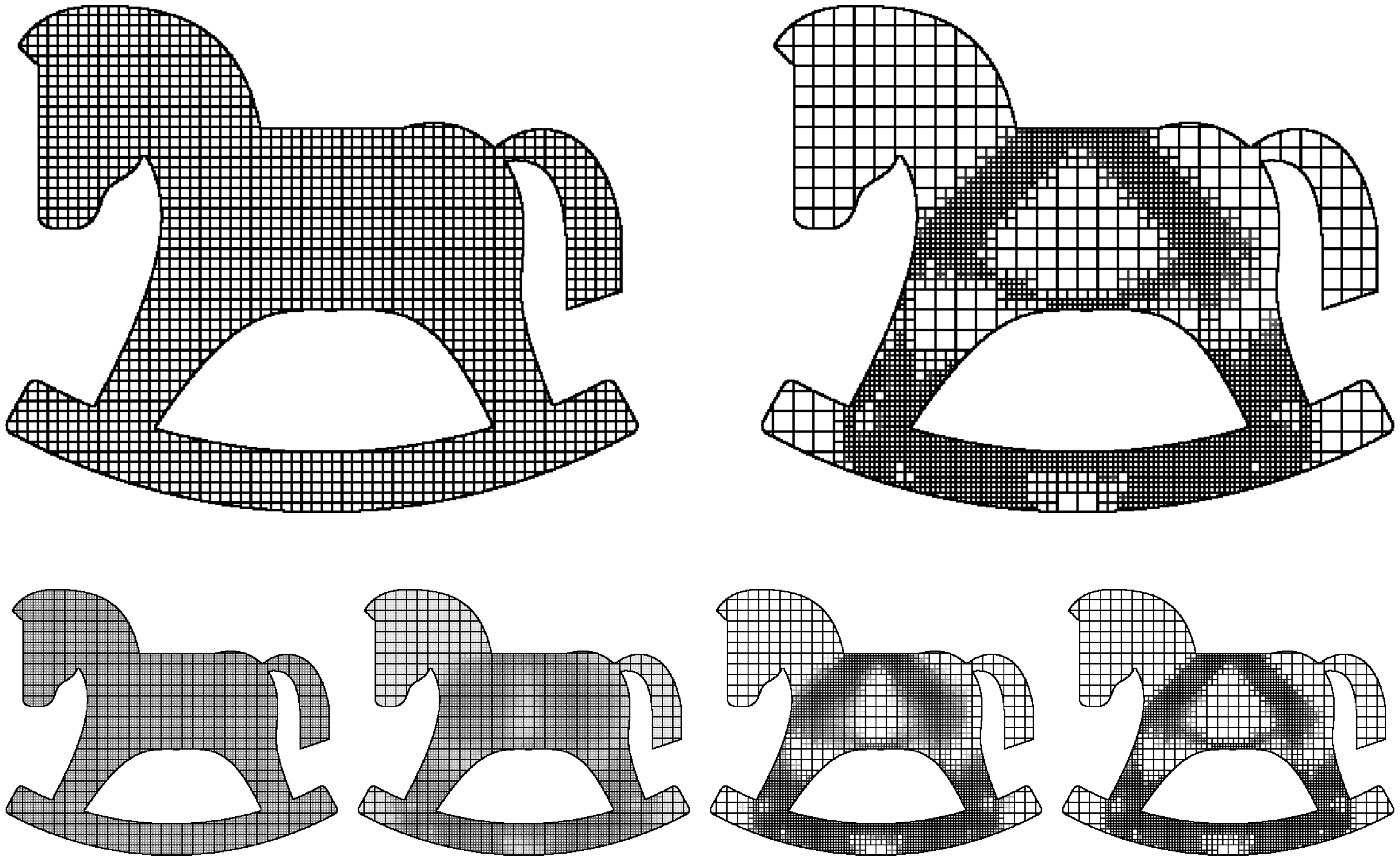
\caption{The rocking horse is filled with a uniform rectangular pattern (a) and the continuously optimized quadtree structure (b). The adaptive quadtree is three times stiffer than the uniform pattern, regarding the prescribed loads which are indicated by blue arrows; The compliances of the uniform and adaptive pattern are $131.6$ and $43.5$, respectively. The two designs use the same amount of material. The bottom row shows the continuous quadtree at four iterations. 
}\label{fig:horse}
\end{figure}

In 3D printing the interior of 3D models is often filled with repetitive geometric patterns. A typical rectangular pattern in a rocking horse model is illustrated in Fig.~\ref{fig:horse} (a). The sparsity of the pattern affects the amount of material usage and the stiffness of the fabricated shape; A denser pattern results in a mechanically stiffer print while consuming more material. To design lightweight and stiff prints, an intuitive option is to start with a sparse pattern and selectively subdivide the cells according to stress analysis, leading to an adaptive structure.

The problem under consideration is the finding of a quadtree structure that maximizes the stiffness regarding prescribed mechanical loads, under the constraint of a given material budget. This is a \textit{discrete} optimization problem: for each cell, to refine or not to refine it. Accurately solving discrete optimization problems is challenging, especially when the number of design variables is large~\cite{Rajeev1992}. Moreover, in (discrete) quadtree optimization the number of design variables is not constant: new cells (and thus design variables) are created as the refinement progresses. A greedy approach to the discrete quadtree optimization problem is suggested in~\cite{Wu16Rhombic} where (rhombic) cells are selectively refined based on a heuristic criterion. While it has been demonstrated that the greedy approach can find a quadtree structure that is stiffer than a uniform pattern with the same amount of material, it is known that heuristic refinements result in a local optimal solution, and it might be away from the global optimum~\cite{Wu16Rhombic}.

To address the limitations of the greedy approach, in this paper we present a novel \textit{continuous} optimization method to the \textit{discrete} problem of quadtree optimization. We relax the binary design variables in optimization to continuous values between $0$ and $1$. By penalizing intermediate values, the continuous optimization eventually converges to a black-white binary design. In the discrete quadtree setting, a cell can only be considered for refinement if it has been created by refining its parent cell. This dependency is encoded in the continuous quadtree setting by a continuous function which smoothly filters out configurations that violate this dependency. Our numerical results demonstrate a good convergence of the continuous formulation. A continuously optimized quadtree is shown in Fig.~\ref{fig:horse}~(b), while a sequence during the optimization is shown in the bottom row. The continuously optimized quadtree structures perform much stiffer than uniform patterns and the heuristically optimized counterparts. In summary, our key contributions are as follows:
\begin{itemize}
\item We extend quadtrees by assigning multi-level continuous design variables to indicate their refinement. This extension, which can be called a \textit{continuous quadtree}, lends itself to gradient-based numerical optimization.
\item We propose a refinement filter to encode the dependency of continuous design variables among multiple levels in the quadtree hierarchy.
\item We demonstrate the effectiveness of quadtree optimization for designing lightweight and stiff infill structures for 3D printing. 
\end{itemize}

The remainder of this paper is organized as follows. After reviewing related work in Section~\ref{sec:RelatedWork}, in Section~\ref{sec:Formulation} we present the continuous quadtree optimization method, and in Section~\ref{sec:Extensions} two extensions. Numerical and physical tests are presented in Section~\ref{sec:Results}, before conclusions are drawn in Section~\ref{sec:Conclusion}.

\section{Related Work}\label{sec:RelatedWork}

We consider quadtree optimization as a \revise{}{structural} topology optimization problem, as the topology \revise{}{of structures} represented by the quadtree changes during the design optimization process. 
\revise{}{Topology optimization aims at finding the material layout which yields for instance the highest stiffness under given constraints. In contrast to shape optimization where the topology is prescribed and remains constant, topology optimization allows for the creation of new voids.} For a thorough review of topology optimization techniques, let us refer to the survey articles~\cite{Sigmund13,Deaton14}. 

Optimizing spatially adaptive structures clearly distinguishes from previous use of adaptive mesh refinement in topology optimization, where the motivation is to reduce the intensive finite element analysis by reducing the number of elements. Spatial adaptive meshes (i.e., quadtree in 2D and  octree in 3D) are often employed in numerical analysis to attain a required accuracy for a minimum amount of computation~\cite{Berger84}. In the context of topology optimization, Maute and Ramm~\cite{Maute95} are among the first to employ adaptive techniques to decrease the number of design variables and to generate smooth structures. Some recent developments along this direction include adaptive polygonal elements (e.g.,~\cite{Nguyen-Xuan2017}), and error control in adaptive meshing (e.g.,~\cite{Wang14ATO,Lambe2018}). The adaptive mesh used in the literature represents the finite elements for elasticity analysis, and its refinement criterion is based on estimated errors. The obtained structures are similar, if not identical, to the structures optimized with uniform fine meshes. In contrast, here the adaptive mesh, its edges in particular, represents the physical structures to be obtained from optimization. The structural refinement respects global measures on the structural stiffness and the material volume. The structures our method is aiming at spread across the entire closed-walled design domain (cf Fig.~\ref{fig:horse} (b)). 

One application of quadtree optimization is to obtain lightweight and stiff infill structures for 3D printing. Topology optimization has been recognized as an important design method for 3D printing, as it fully exploits the manufacturing flexibility to create lightweight and mechanically optimal structures~\cite{Brackett2011}. A recent focus in this area has been on incorporating manufacturability constraints (e.g., overhang angle~\cite{Langelaar16,Gaynor16,Mirzendehdel2016CAD,Qian2017} and length scale~\cite{Lazarov16}) and thus eliminating post-processing of the numerically optimized structures. 

Infill is not a constraint in 3D printing, but rather a feature that can be exploited for improving mechanical stability~\cite{Clausen16,Wu17-infill}. A summary of inner structures and their optimization is presented in the survey paper by Livesu et al.~\cite{Livesu2017}. The mechanically optimized structures take different forms, including frames~\cite{Wang13,Zhang2015,Jiang2017}, honeycomb-like structures~\cite{Lu14}, non-uniform shells~\cite{Musialski2016,Zhao2017,Li2017CGF}, foams~\cite{Martinez16,Martinez2017}, micro- and lattice  structures~\cite{Schumacher15,Panetta15,Chougrani2017,Zhu2017} and bone-mimicking structures~\cite{Wu17-infill,Wu2017CMAME,Liu2017}. In this work we focus on subdivision structures which allow for an intuitive control over some geometric features. For instance, the interior quadtree structures have a uniform thickness, which simplifies tool-path generation in 3D printing. The cells have a fixed and prescribed orientation, which has been used to ensure self-supportability~\cite{Wu16Rhombic}. Furthermore, the coarsest and finest grids effectively determine the maximum and minimum void sizes, respectively. The uniform thickness and the minimum void size have positive implications on thermal effects in the additive manufacturing process~\cite{Gibson10}. Going beyond 3D printing, structures with a certain level of regularity facilitate the production of large engineering structures in architectural geometry~\cite{Pottmann2007}. In this case, the final quadtree can be interpreted as being assembled from frames. The frames have only a few variations in size, and thus each variation can be massively produced~\cite{Norato2015}.

\section{Continuous Quadtree Optimization}\label{sec:Formulation}

Given a 2D design domain, our method aims at finding the quadtree structure that is stiffest regarding prescribed mechanical loads on the domain boundary. \revise{}{The design domain is defined by an input shape. In this section we take a simple rectangular shape to explain the algorithm.} The optimization is subject to a constraint on the limit of material usage. Before we formulate this optimization problem in Section~\ref{subsec:Optimization}, we present in Section~\ref{subsec:Discretization} the design parameterization, including the design variables of quadtree refinements, the finite elements for structural analysis, and the mapping from design variables to finite elements. A refinement filter encoding the dependency of design variables is proposed in Section~\ref{subsec:RefinementFilter}.

\subsection{Design Parameterization}\label{subsec:Discretization}

\begin{figure}[t]
\centering
\small
\def\svgwidth{0.64\linewidth}
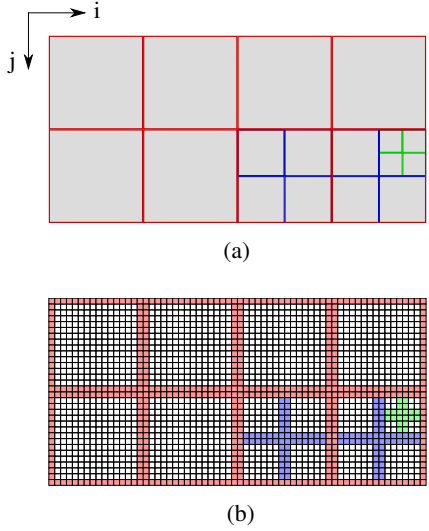
\caption{Illustration of design parameterization \revise{}{on a rectangular shape} using the quadtree grid (a) and the underlying finite elements (b) which are used for evaluating mechanical response of the design.
}\label{fig:grid}
\end{figure}

We formulate quadtree optimization as a continuous material distribution problem. To this end, a hybrid discretization is employed. As illustrated in Fig.~\ref{fig:grid} for a rectangular design domain, the hybrid discretization is composed of a quadtree grid which encodes the topology information, and a discretization by uniform square elements which encode material distribution and upon which elasticity analysis is performed. In Fig.~\ref{fig:grid}~(a) the red lines form the initial coarse grid where refinements start from. The coarse cells are mapped to blocks of square elements in the finite element grid in Fig.~\ref{fig:grid}~(b). 
Within each block, elements in the first boundary layer are assigned solid (i.e., a density value of $1.0$), indicated by the red-shaded elements. Adjacent coarse cells share a common edge, and this edge in the finite element grid has a thickness of two elements --- one element on each side. 

For illustrative purpose a few cells in the bottom right of the domain are refined, creating \revise{cross-shaped}{plus-shaped (''+'')} sub-structures. The \revise{crosses}{plus-shapes} are also mapped to solid elements, with colours indicating the different levels of refinement. 
In this example two refinement levels are applied. A further refinement level would merge the edges and form fully solid blocks in the finite element grid. While such fully solid blocks cause no problem for the optimization, here we exclude this situation for the purpose of ensuring a uniform thickness inside the design domain. 

\subsubsection{Multi-level Design Variables ($x^k$)}
To enable numerical optimization, we extend quadtree by assigning a continuous design variable $x_{i,j}\in[0,1]$ to each of the cells in the quadtree grid, with $x_{i,j}=1$ (resp. $x_{i,j}=0$) referring to a refinement (resp. non-refinement) for the cell with the 2D indices $(i,j)$ \revise{}{in the 2D grid (Fig.~\ref{fig:grid})}. The design variable ($x$) effectively indicates the intensity of the \revise{cross-shaped}{plus-shaped} sub-structure that will be created if the cell under consideration is refined. Design variables are assigned on each level in the quadtree, denoted by $x^k_{i,j}$, where $k$ indicates the level. The initial coarse grid is denoted as the first level. From the first level downwards, the resolution of cells, which is also the resolution of the design field, is doubled on each dimension from the previous level. Denoting the resolution on the $k$-th level by $(n^k_x, n^k_y)$, and the resolution of the initial coarse grid by $(n^0_x, n^0_y)$, the relation can be written as
\begin{equation}
(n^k_x, \; n^k_y) = (2^{k-1}\,n^0_x, \; 2^{k-1}\,n^0_y), \qquad \text{$k=1,...,\bar{k}$},
\label{eq:Resolution}
\end{equation}
where $\bar{k}$ is the maximum allowed refinement level. For a coarse block of $2^m \times 2^m$ elements, $\bar{k}$ is calculated by
\begin{equation}
\bar{k} = m-2.
\end{equation}
When the refinement level reaches $\bar{k}$, it creates voids in the size of $2^2$ finite elements.

\subsubsection{Mapping ($x^k \to \rho$)}
The mapping from the design variables ($x^k$) on the quadtree grid to the density field ($\rho$) in finite element analysis has been illustrated in Fig.~\ref{fig:grid}. The density of finite elements is mapped from the design variables $x^k$, $k=1,...,\bar{k}$, corresponding to \revise{cross-shaped}{plus-shaped} sub-structures (blue and green in Fig.~\ref{fig:grid}). The boundary of the initial coarse cells also maps to density of finite elements (red in Fig.~\ref{fig:grid}). For compactness in the following formulation, let us denote the initial coarse cells by $x^0$ which has the same resolution as the design variable on the first level $x^1$, i.e., $(n^0_x,n^0_y) = (n^1_x,n^1_y)$. Corresponding to a solid coarse frame for ensuring a minimum stiffness of the quadtree structure, $x^0$ has a constant value $x^0_e = 1$ for all coarse cells.

Reshaping the multi-level design fields $x^k|_{n^k_x \times n^k_y}$, $k=0,...,\bar{k}$, and the density field $\rho|_{2^m n^0_x \times 2^m n^0_y}$ into column vectors $\bm{x}^k$ and $\bm{\rho}$, respectively, the mapping can be written compactly by a transformation,
\begin{equation}
\bm{\rho} = \textstyle{\sum^{\bar{k}}\limits_{k=0}} \bm{T}^k \bm{x}^k.
\label{eq:Mapping}
\end{equation}
The transformation $\bm{T}^k$ is a sparse matrix with a dimension of $2^{2m} n^0_x n^0_y \times n^k_x n^k_y$. The non-zero entries in $\bm{T}^k$ have a constant value of $1$, and concern the specific finite elements that will be assigned by refinements $x^k_{i,j}$. Since the non-zero entries have a unit value, the transformation assigns the same value of the refinement to the affected finite elements. For instance, an intermediate refinement by $x^k_{i,j} = 0.5$ leads to intermediate densities for the corresponding finite elements, $\rho=0.5$.

\subsection{Optimization Problem}\label{subsec:Optimization}

With the multi-level continuous refinement field ($x^k$) and the physical density field ($\rho$) defined in the previous subsection, we can write the compliance minimization problem as follows: 
\begin{align}
\underset{\bm{x}}{\text{min}}
&\quad c = \bm{U}^T \bm{K}(\bm{\rho}) \bm{U}, \label{eq:objective} \\
\text{s.t.}
&\quad \bm{K}(\bm{\rho}) \bm{U} = \bm{F},  \label{eq:state} \\
&\quad V(\bm{\rho}) = \textstyle{\sum\limits_{\forall e}} \rho_e v_e \le V^*, \label{eq:globalConstraint} \\
&\quad x^k_{i,j} \in [0,1]. \label{eq:design}
\end{align}
In contrast to classic topology optimization~\cite{Sigmund01,Wu16} where the design variable is the density field ($\rho$), for obtaining quadtree structures we consider the refinement field ($x^k$) as new design variables. The density vector $\bm{\rho}$ is calculated from $\bm{x}^k$ by Eq.~\eqref{eq:Mapping}.

The objective is to minimize the compliance ($c$), which is equivalent for stiffness maximization. The first constraint is the static equation of linear elasticity. $\bm{K}$, $\bm{U}$, and $\bm{F}$ are stiffness matrix, displacement vector, and force vector, respectively. The global stiffness matrix $\bm{K}$ is assembled from element stiffness matrix $\bm{k}_e$. The modified SIMP interpolation model~\cite{Andreassen10} is used to compute the element stiffness matrix: 
\begin{align}
\bm{k}_e &= E_e(\rho_e) \bm{k}_0, \label{eq:elementStiffness} \\
E_e(\rho_e) &= E_{\text{min}} + \rho^{p}_e \, (E_0 - E_{\text{min}}), \label{eq:YoungsModulus}
\end{align}
where $\bm{k}_0$ is the element stiffness matrix for an element with unit Young's modulus. $E_0$ is the stiffness of the material. $E_{\text{min}}$ is a very small stiffness assigned to void elements in order to prevent the stiffness matrix from becoming singular. And $p$ is a penalization factor (typically $p = 3$) which is introduced to ensure black-and-white solutions~\cite{Sigmund01}. The second constraint measures the volume occupied by solid elements. $v_e$ is the unit volume for each element, and $V^*$ is the maximum allowed volume. The last constraint restricts the design variables to take a value between $0$ (non-refinement) and $1$ (refinement). 

\begin{figure}[t]
\centering
\includegraphics[width=0.98\linewidth]{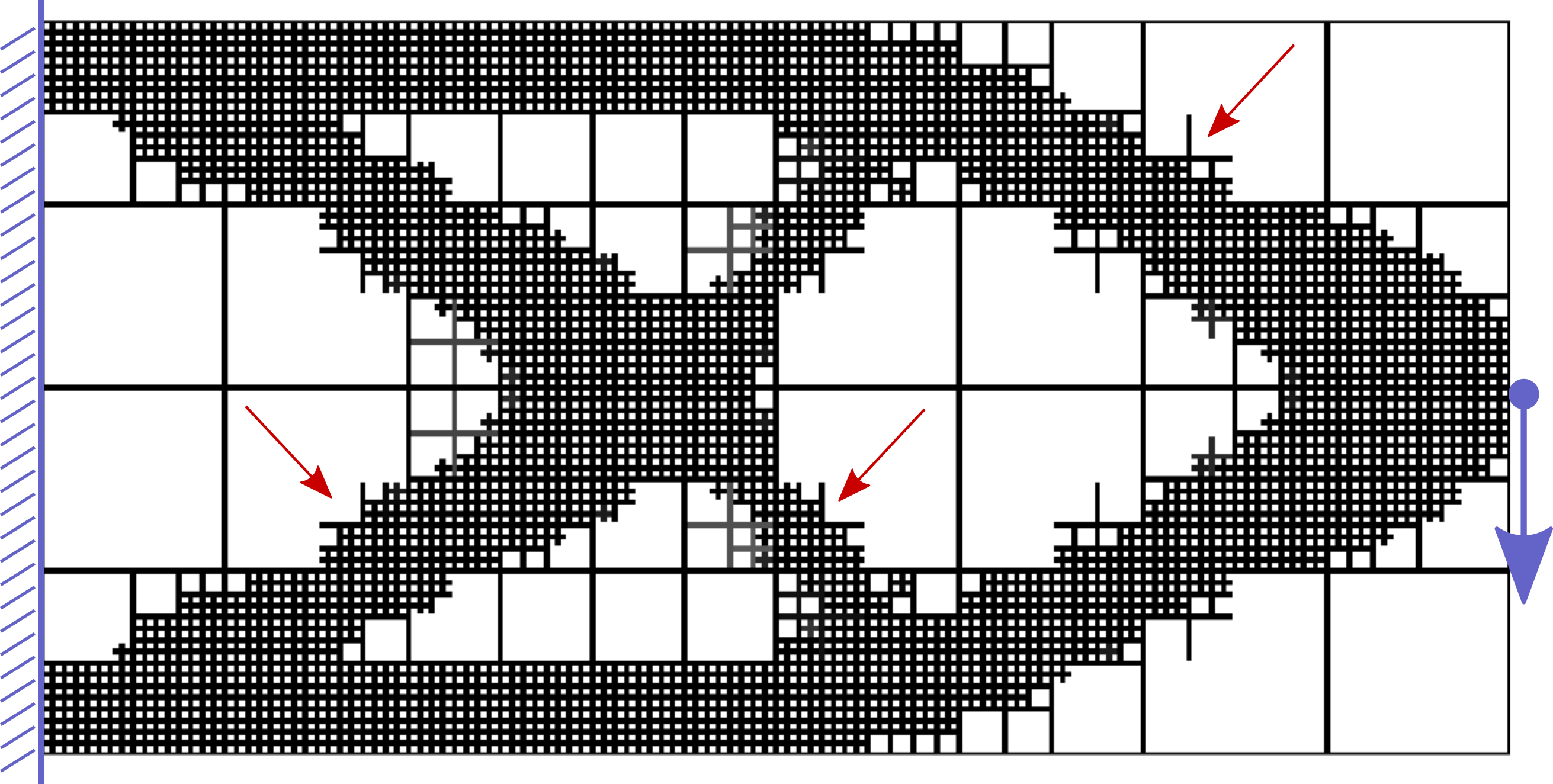}
\caption{When the refinement rules are not applied, the optimized structure consists of many suspended structures, indicated by the red arrows. The compliance value is $136.1$.
}\label{fig:beamNoFilter}
\end{figure}

Solving this optimization problem for a cantilever beam presents the design shown in Fig.~\ref{fig:beamNoFilter}. The beam is fixed on the left side, while a force is applied on the right side in the middle. The prescribed material volume is $40\%$ of the design domain. The optimized structure shows crossing sub-structures which correspond to the design variables. However, many of such sub-structures (pointed out by the red arrows) are suspended. This is due to the fact that the (imaginary) fine cells are refined, even though the fine cells themselves have not yet been created from refining their parent cells. The suspended sub-structures don't effectively carry mechanical loads, and disobey the intention to create quadtree-like mechanical structures. In the following, a refinement filter is proposed to encode the recursive refinement rules, and thus to effectively filter out such suspended sub-structures.

\subsection{Refinement Filter ($x^k \to \tilde{x}^k$)}
\label{subsec:RefinementFilter}

In creating a quadtree the refinement is only applicable to a cell that already exists. This implies dependency of the refinement of a fine cell on the refinement of its parent cell. In the binary setting, the refinement of a cell on the $k-$th level, $x^k_{i,j}$, will not happen, as long as the refinement of its parent cell has not yet taken place, i.e., $x^{k-1}_{i^{-1},j^{-1}}=0$. Here $(i^{-1},j^{-1})$ refers to the indices of the parent cell. To encode this dependency in the setting of continuous design variables, we introduce a filtered design variable $\tilde{x}^k_{i,j}$:
\begin{equation}
\tilde{x}^k_{i,j} = \text{min}{(x^k_{i,j}, \tilde{x}^{k-1}_{i^{-1},j^{-1}})}.
\label{eq:Min}
\end{equation}
In case of $x^{k-1}_{i^{-1},j^{-1}} = 0$, i.e., the parent cell is not refined, the filtered value is
\begin{equation}
\tilde{x}^k_{i,j} = \text{min}{(x^k_{i,j}, 0)} = 0.
\end{equation}
This $\text{min}$ function thus restricts the design variable on the fine level when the parent cell has not been refined. In case of $x^{k-1}_{i^{-1},j^{-1}} = 1$, the filtered value is unaffected,
\begin{equation}
\tilde{x}^k_{i,j} = \text{min}{(x^k_{i,j}, 1)} = x^k_{i,j}.
\end{equation}

\begin{figure}[t]
\centering
\small
\def\svgwidth{0.9\linewidth}
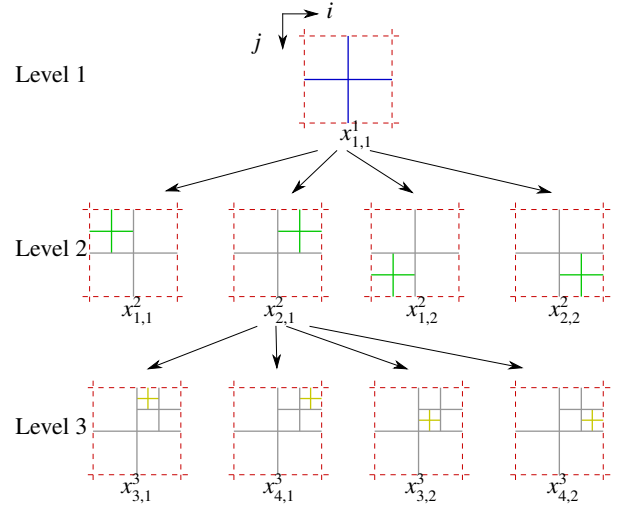
\caption{Illustration of dependency in recursive refinement. For instance, $x^3_{4,1}$ depends on $x^2_{2,1}$ which depends on $x^1_{1,1}$.
}\label{fig:refinement}
\end{figure}

Applying refinement recursively creates a quadtree with multiple levels. The dependency thus goes from the fine cell through all the intermediate levels to one of the coarse cells on the first level. For instance, consider the design variable $x^3_{4,1}$ illustrated in Fig.~\ref{fig:refinement}. The filtered value is 
$\tilde{x}^3_{4,1}=\text{min}{(x^3_{4,1}, \tilde{x}^2_{2,1})}$, where $\tilde{x}^2_{2,1}=\text{min}{(x^2_{2,1}, \tilde{x}^1_{1,1})}$. At first glance this requests multiple filters executed sequentially, leading to much complex sensitive analysis afterwards. However, the recursive $\text{min}$ functions can be consolidated into an equivalent $\text{min}$ function by taking in as input all the involved arguments in the individual $\text{min}$ functions, i.e.,
\begin{equation}
\tilde{x}^k_{i,j} = \text{min}{\left(x^k_{i,j}, \, x^{k-1}_{i^{-1},j^{-1}}, \,..., \, x^{1}_{i^{-k+1},j^{-k+1}} \right)}.
\label{eq:minExpand}
\end{equation}
This consolidation does not introduce approximate error, and significantly simplifies sensitivity analysis (cf. Appendix).

After the design variables have been filtered, the density vector, $\bm{\rho}$, is updated. This is done by replacing design vector $\bm{x}^k$ by the filtered design vector $\tilde{\bm{x}}^k$ in Eq.~\eqref{eq:Mapping}, i.e.,
\begin{equation}
\bm{\rho} = \textstyle{\sum^{\bar{k}}\limits_{k=0}} \bm{T}^k \tilde{\bm{x}}^k.
\label{eq:MappingFilter}
\end{equation}

\begin{figure*}[t]
\centering
\includegraphics[width=0.98\linewidth]{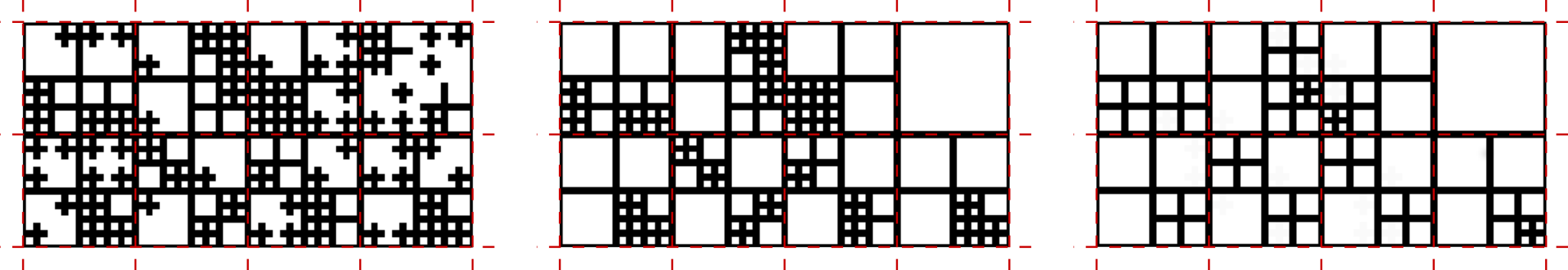}
\caption{Left: A randomly initialized design field shows many suspended structures. Middle: The suspended structures disappear when the refinement filter is applied. Right: A balanced quadtree obtained by filtering out fine cells which are at least two levels smaller than one of its neighbours. 
}\label{fig:refinementFilter}
\end{figure*}

\begin{figure}[t]
\centering
\includegraphics[width=0.98\linewidth]{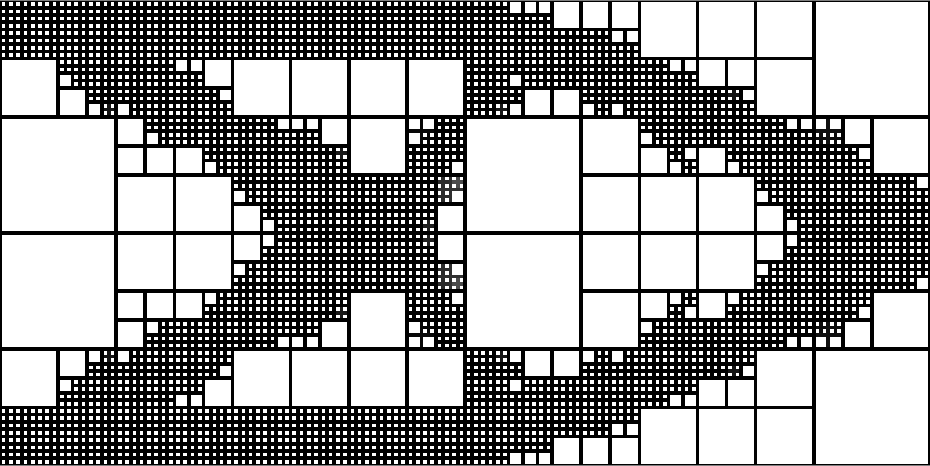}
\caption{With the refinement filter the optimization creates a quadtree structure.  
The compliance value is $143.3$.
}\label{fig:quadtreeBeamUnbalanced}
\end{figure}

The effect of the recursive refinement filter is visualized in Fig.~\ref{fig:refinementFilter}. On the left, the multi-level design variables are initialized randomly with a value of $0$ or $1$. A few suspended sub-structures can be observed in the left sub-figure, and they disappear in the middle sub-figure. Fig.~\ref{fig:quadtreeBeamUnbalanced} shows the result of integrating this filter into the optimization for the cantilever beam. 
It is worth noting that Fig.~\ref{fig:quadtreeBeamUnbalanced} is not post-processed from Fig.~\ref{fig:beamNoFilter}; The filter acts as a constraint that is implicitly integrated into the iterative optimization process. 

\subsubsection{Smooth Approximation}
Gradient-based numerical optimization necessitates differentiable functions. To this end, the non-differentiable $\text{min}$ function (Eq.~\eqref{eq:minExpand}) is approximated by a smooth $p$-norm function, with a negative exponent,
\begin{align}
\tilde{x}^k_{i,j} &= \text{min}{\left(x^k_{i,j}, \, x^{k-1}_{i^{-1},j^{-1}}, \,..., \, x^{1}_{i^{-k+1},j^{-k+1}} \right)} \\
&\approx \left( \frac{1}{k} \, \textstyle{\sum^{k-1}\limits_{l=0}}  \left(x^{k-l}_{i^{-l},j^{-l}} \right)^{p_n} \right)^{\frac{1}{p_n}}.
\label{eq:filter}
\end{align}
Here $p_n$ represents the $p$-norm exponent, to be distinguished from the penalty $p$ in SIMP interpolation (Eq.~\eqref{eq:YoungsModulus}). As $p_n$ approaches negative infinity, the $p$-norm equals the minimum value.
The $p$-norm is normalized by a factor considering the number of arguments, $\frac{1}{k}$ in this formulation. The normalization reduces the approximation error when a practical $p_n$ value is not infinite~\cite{Wu17-infill}. In our tests, we use $p_n=-16$.

The indices of parent cells are needed in Eq.~\eqref{eq:filter}. When using a regular discretization, the indices can be found by a recursive function,
\begin{equation}
i^{-l} = h(i,-l) = h(h(i,-l+1),-1),
\end{equation}
where the function $h$ is defined as
\begin{equation}
h(i,-1) = \text{floor}\left(\frac{1}{2}\left(i+1 \right)\right).
\end{equation}

\section{Extensions}\label{sec:Extensions}

\subsection{Balanced Quadtree}\label{subsec:Balanced}
We provide an option to restrict the quadtree such that the level difference between neighbouring cells is at most one, known as balanced or restricted quadtrees. A balanced quadtree maintains a smooth transition from fine cells to coarse cells. For 3D printed infill, a smooth transition is expected to make the part more stable for uncertain loads. 

In an unbalanced quadtree, the refinement of a cell depends only on its parent cell (Eq.~\eqref{eq:Min}). In the balanced quadtree, the refinement additionally depends on the neighbours of its parent cell. The refinement filter (Eq.~\eqref{eq:Min}) is thus updated by
\begin{equation}
\tilde{x}^k_{i,j} = \text{min}(x^k_{i,j}, \, \tilde{x}^{k-1}_{i^{-1},j^{-1}}, \, \tilde{x}^{k-1}_{i^{-1} \pm 1,j^{-1}\pm1}).
\label{eq:MinBalanced}
\end{equation}
Here, it is assumed that the indices $i^{-1} \pm 1$ belong to $[1,n^{k-1}_x]$ and $ j^{-1}\pm1$ belong to $[1,n^{k-1}_y]$, otherwise the invalid arguments are excluded. A consolidation of this recursive function is obtained similar to Eq.~\eqref{eq:minExpand}. The recursive function can point to the same coarser cell multiple times, due to the introduction of neighbours. Only one copy of each dependent cell is included in the consolidated filter. Since the encoding of neighbourhood and hierarchy remains constant during the optimization, finding the dependency is only performed once in the initialization.

\begin{figure}[t]
\centering
\includegraphics[width=0.98\linewidth]{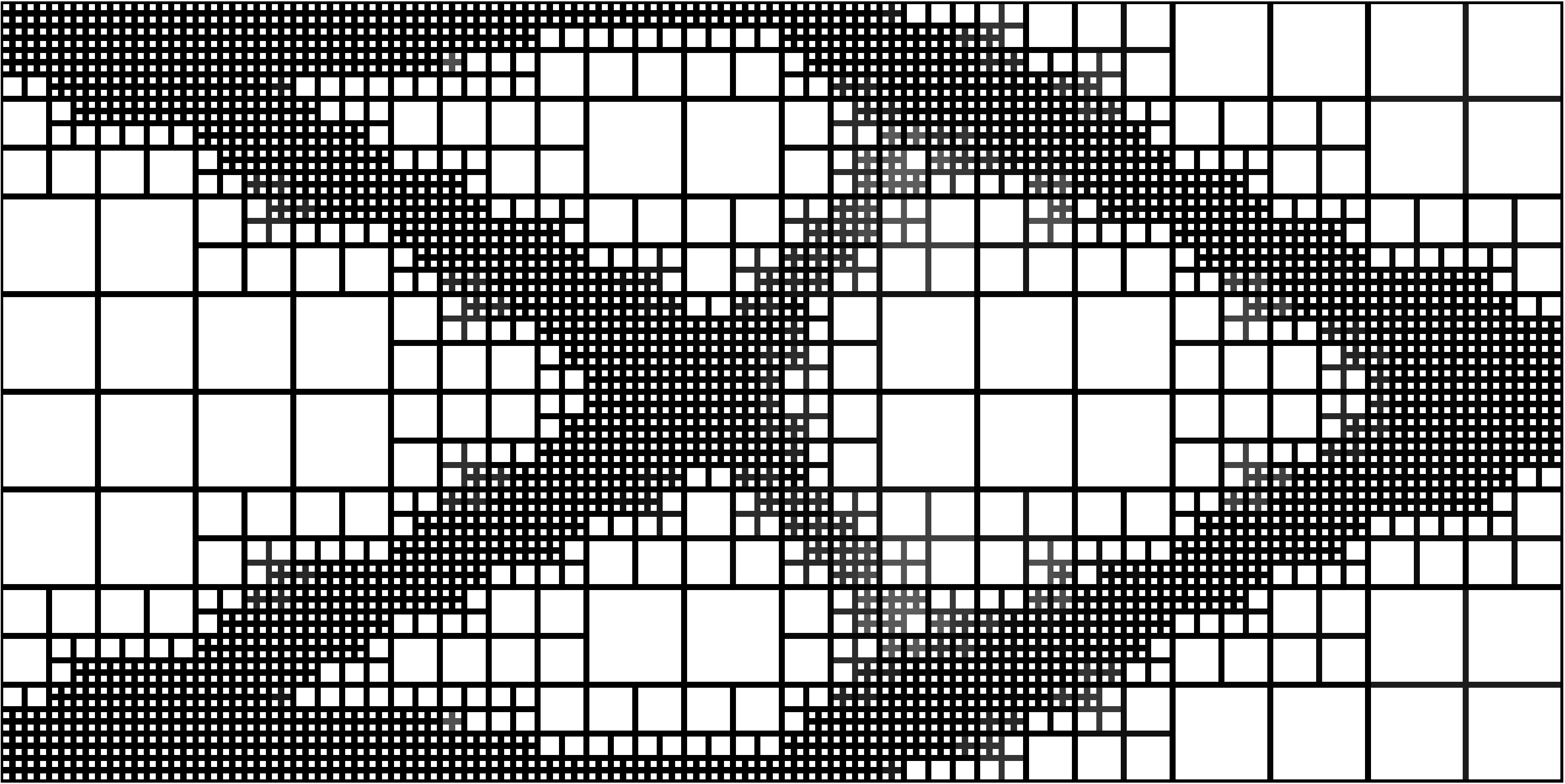}
\caption{A \textit{balanced} quadtree structure obtained from optimizing a cantilever beam. The compliance value is $165.7$.
}\label{fig:balancedBeam}
\end{figure}

The effect of the balanced refinement filter can be observed in Fig.~\ref{fig:refinementFilter} (right). An optimized balanced quadtree for the cantilever beam is shown in Fig.~\ref{fig:balancedBeam}. The three versions (i.e., Fig.~\ref{fig:beamNoFilter}, Fig.~\ref{fig:quadtreeBeamUnbalanced} and Fig.~\ref{fig:balancedBeam}) use the same amount of material. As more restrictions on refinement are introduced, the compliance value becomes larger, meaning that the stiffness is decreased. This moderate sacrifice in stiffness is actually rewarded by a significant increase in robustness, which will be discussed in results section (cf. Fig.~\ref{fig:robust}).

\subsection{Integration with other Filters}\label{subsec:Heaviside}
In density-based topology optimization using square elements, regulations such as density filter~\cite{Guest04} or sensitivity filter~\cite{Sigmund01} are often applied to get rid of checkerboard patterns (i.e., regions of alternating black and white elements). In our formulation, the refinement itself can be regarded as a regulation, and the density/sensitivity filter is no longer necessary.

In topology optimization with continuous design variables, projection filters are often applied to improve the convergence of approaching a binary black-white design. Our continuous reformulation of refinement also benefits from such a projection. 
In particular, we integrate the projection filter proposed in~\cite{Wang10}:
\begin{equation}
\bar{x}^k_{i,j} = \frac{\text{tanh} \left(\beta \eta \right) + \text{tanh} \left(\beta\left({x}^k_{i,j} - \eta \right) \right)}{\text{tanh} \left(\beta \eta \right) + \text{tanh} \left(\beta\left(1 - \eta \right) \right)}.
\label{eq:Projection}
\end{equation}
Here $\beta$ controls the sharpness of projection. In the limit of ${\beta\to\infty}$ the projection approaches a discontinuous step function. We employ a continuation of this parameter. $\eta \in [0,1]$ is the threshold value, and we choose $\eta=0.5$. 

\section{Results}\label{sec:Results}

\begin{figure*}[!htb]
\centering
\subfloat[It. $40$, $c=317.2$, $s=0.631$]
{\includegraphics[width=0.32\linewidth]{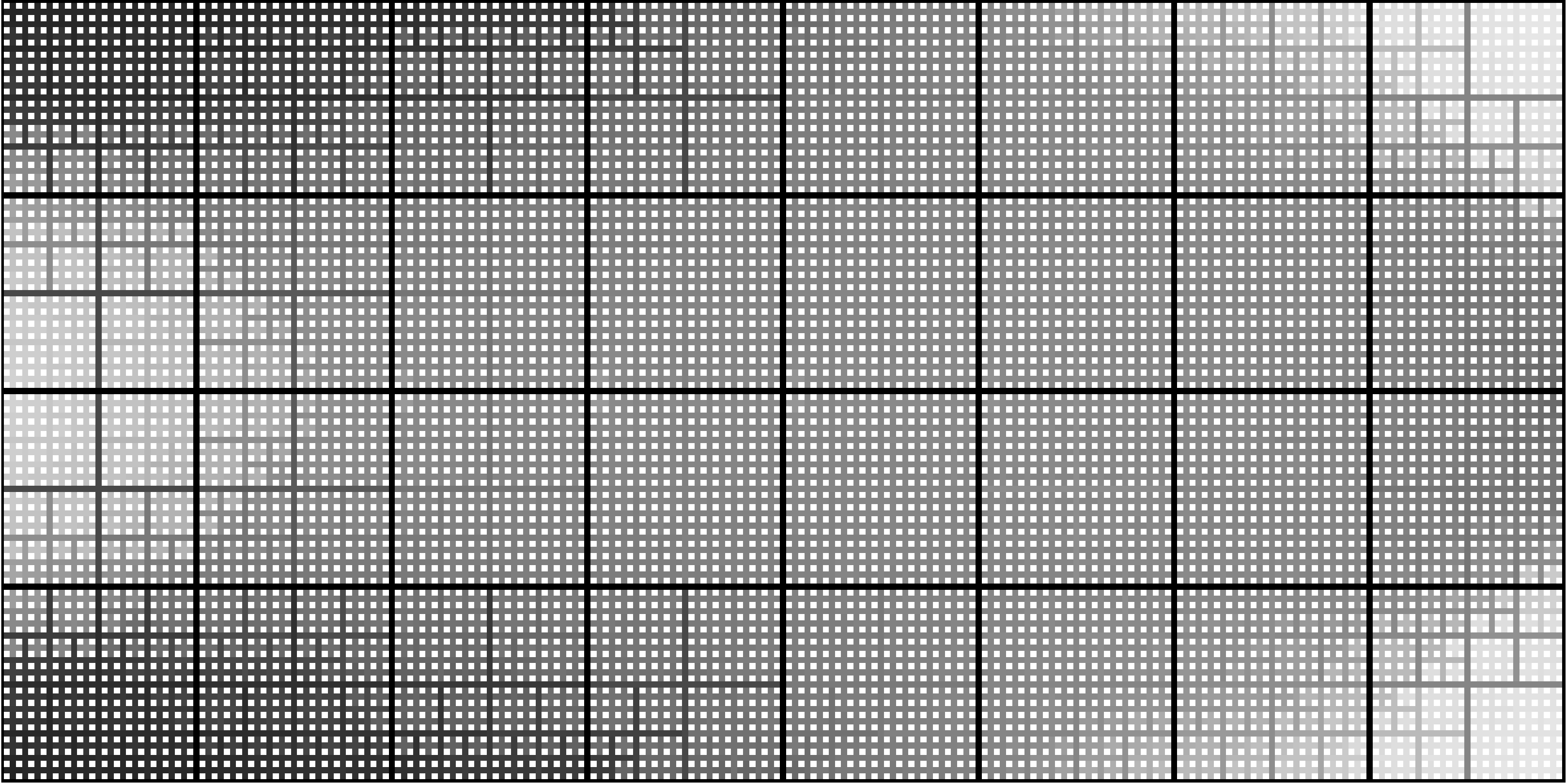}}
\quad
\subfloat[It. $100$, $c=249.8$, $s=0.500$]
{\includegraphics[width=0.32\linewidth]{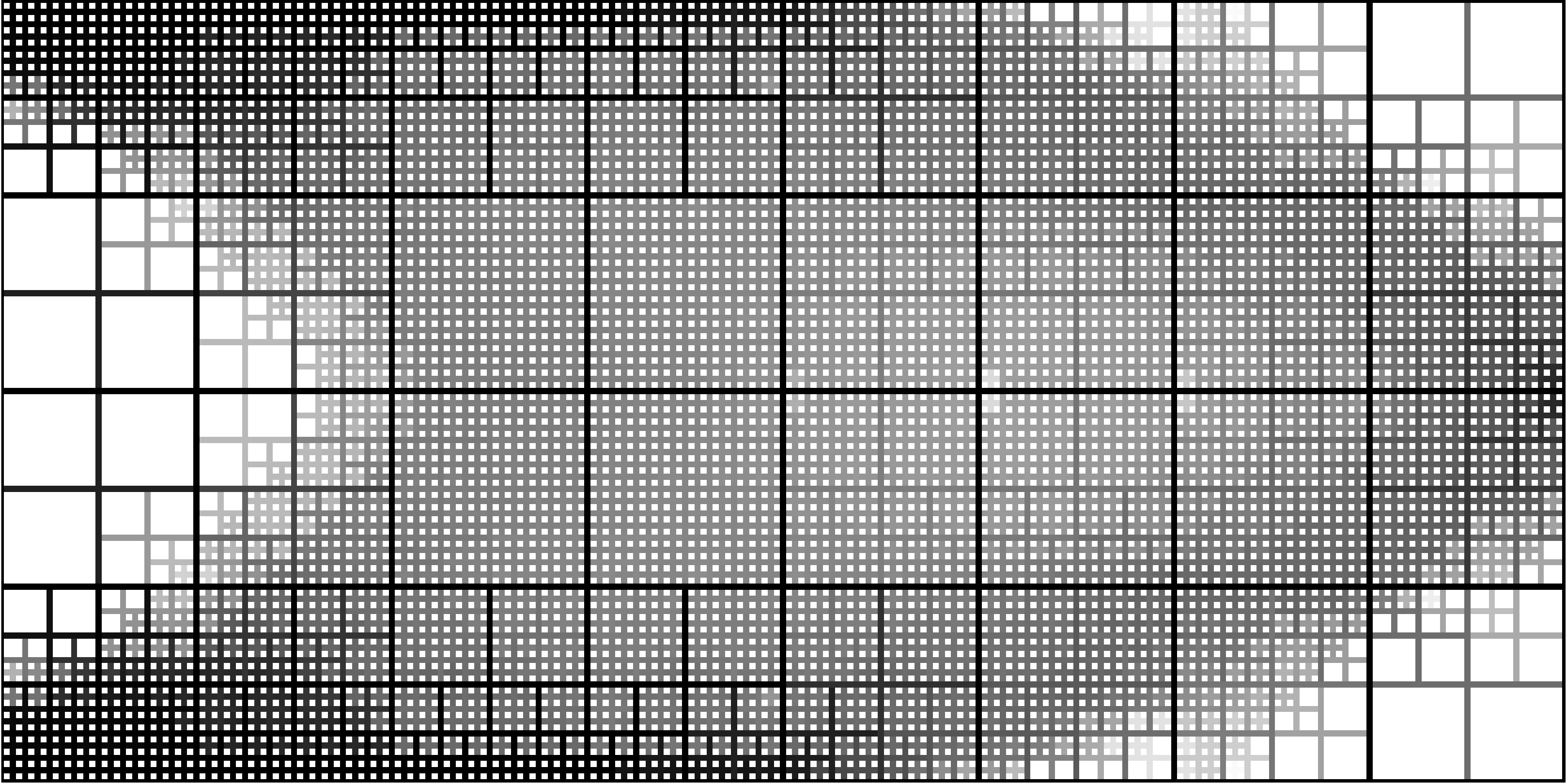}} 
\quad
\subfloat[It. $160$, $c=215.4$, $s=0.362$]
{\includegraphics[width=0.32\linewidth]{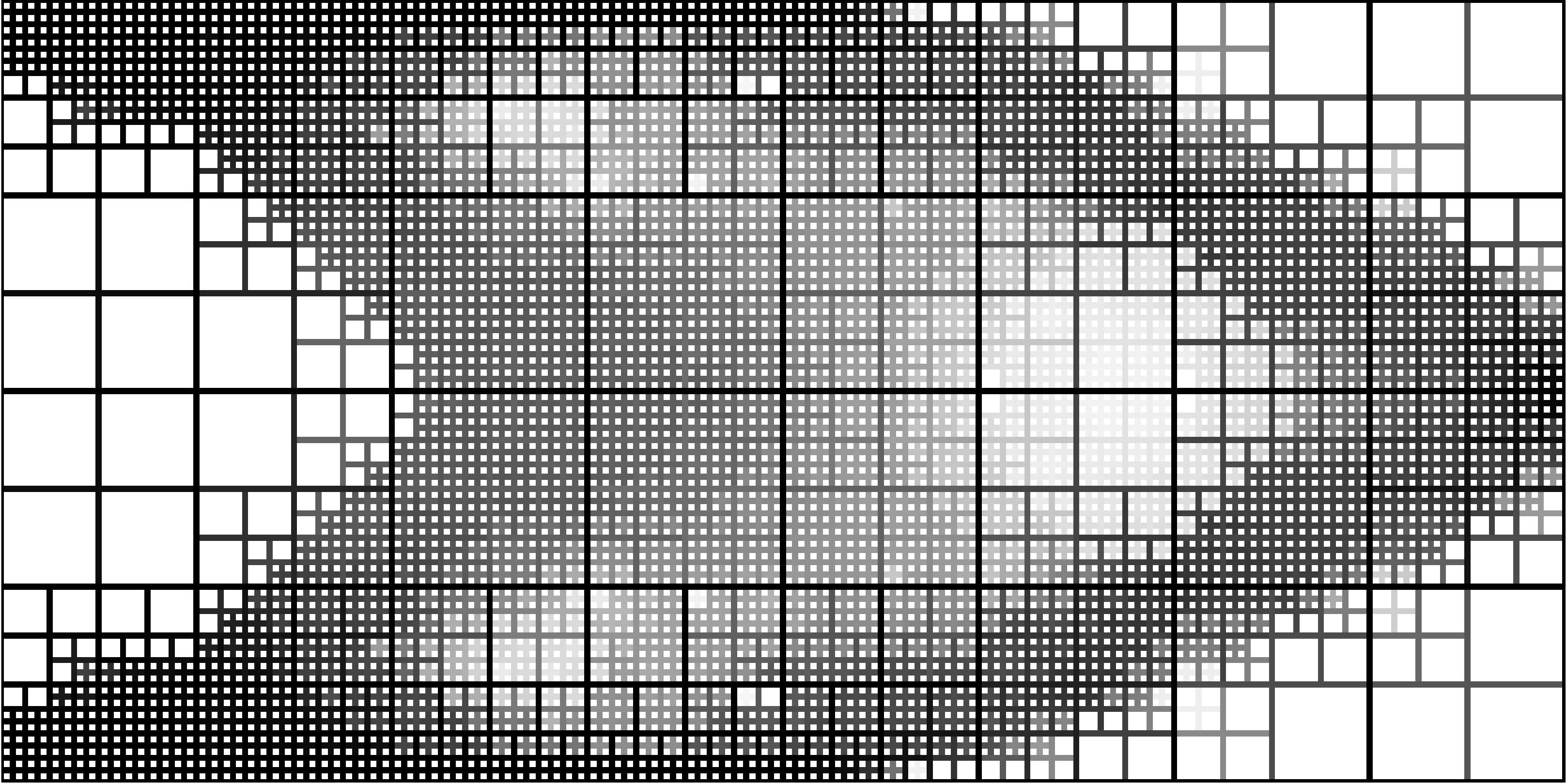}} \\
\subfloat[It. $220$, $c=177.9$, $s=0.151$]
{\includegraphics[width=0.32\linewidth]{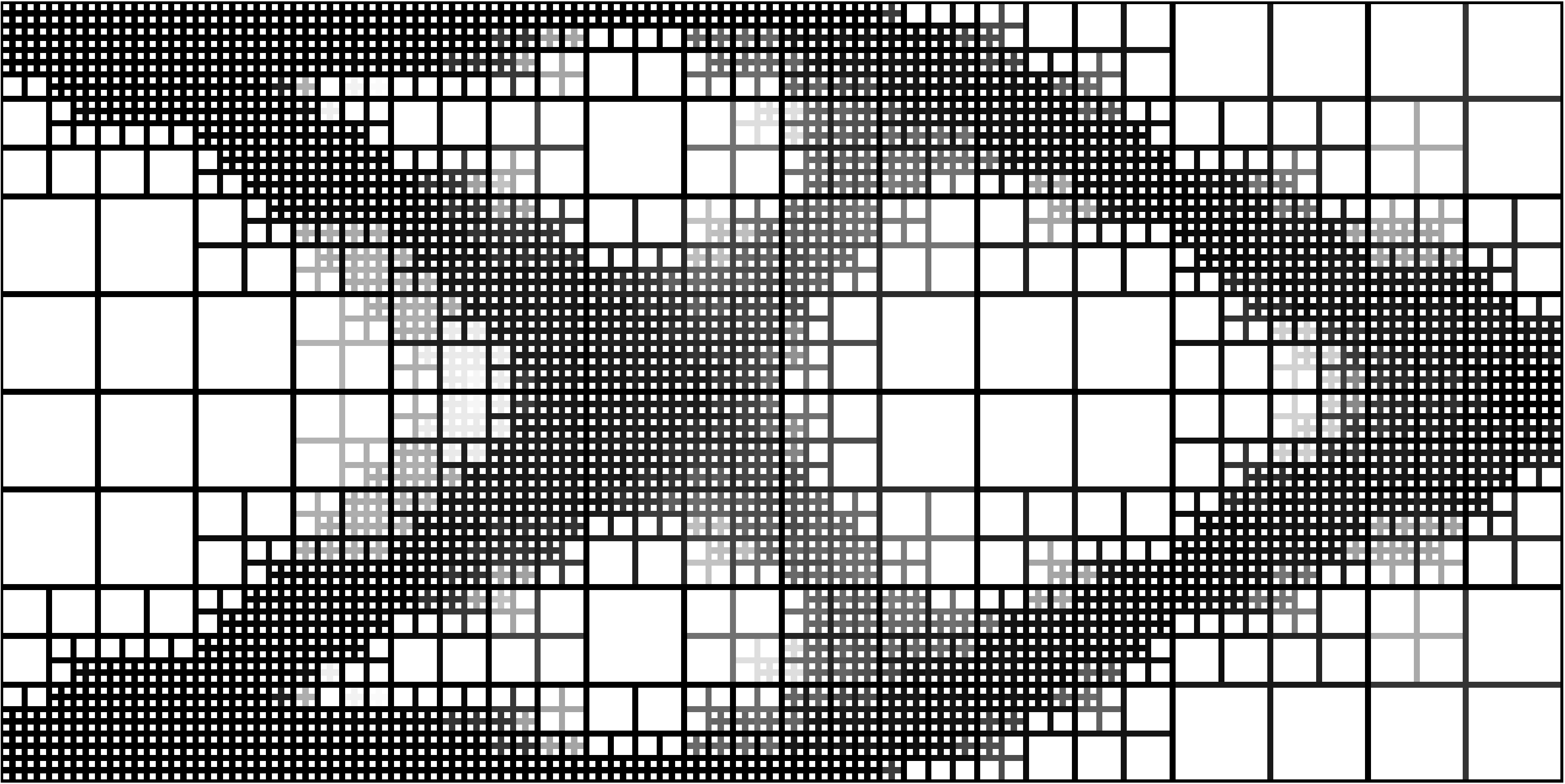}} 
\quad
\subfloat[It. $280$, $c=169.2$, $s=0.070$]
{\includegraphics[width=0.32\linewidth]{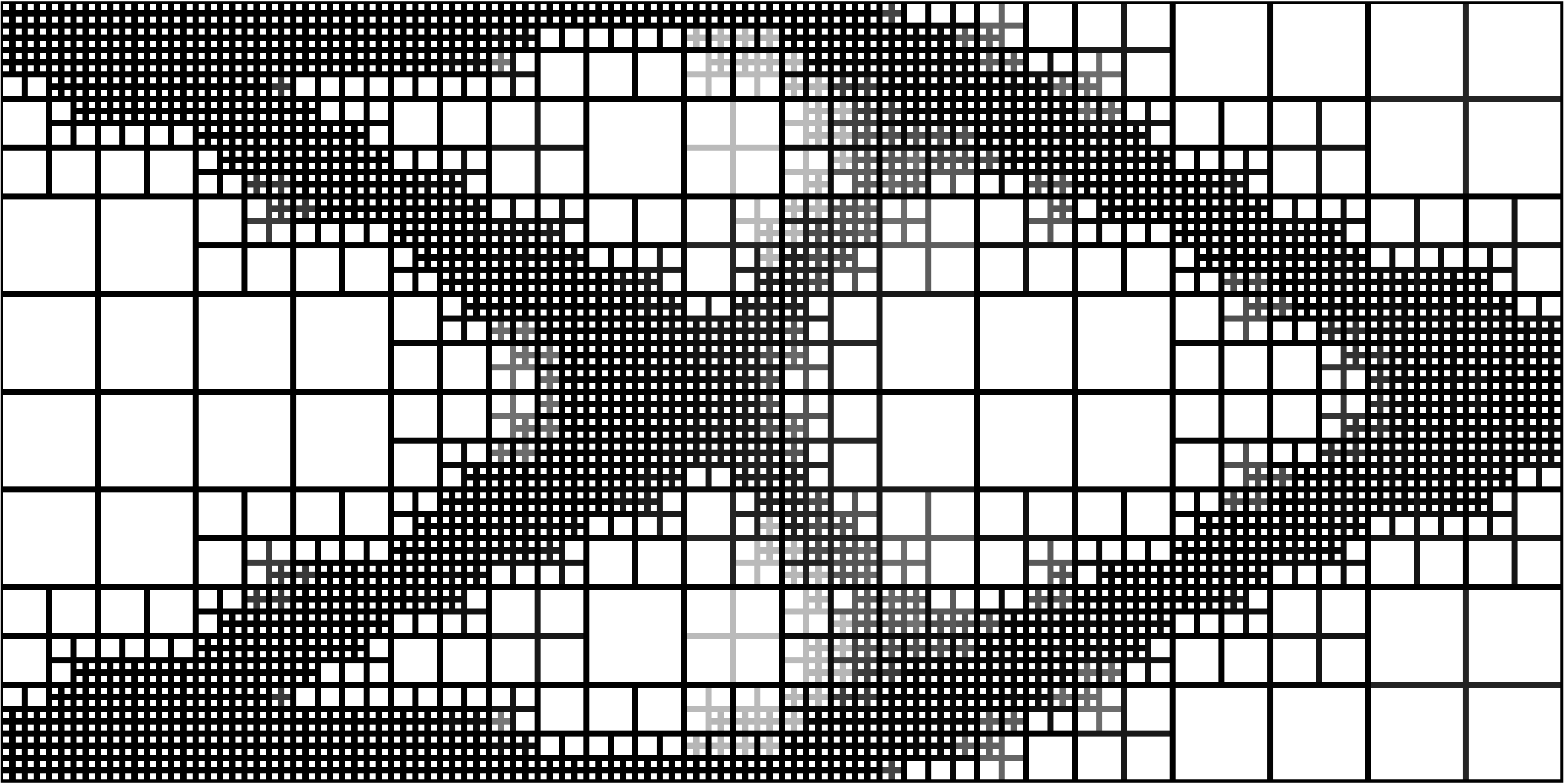}}
\quad
\subfloat[It. $360$, $c=165.7$, $s=0.040$]
{\includegraphics[width=0.32\linewidth]{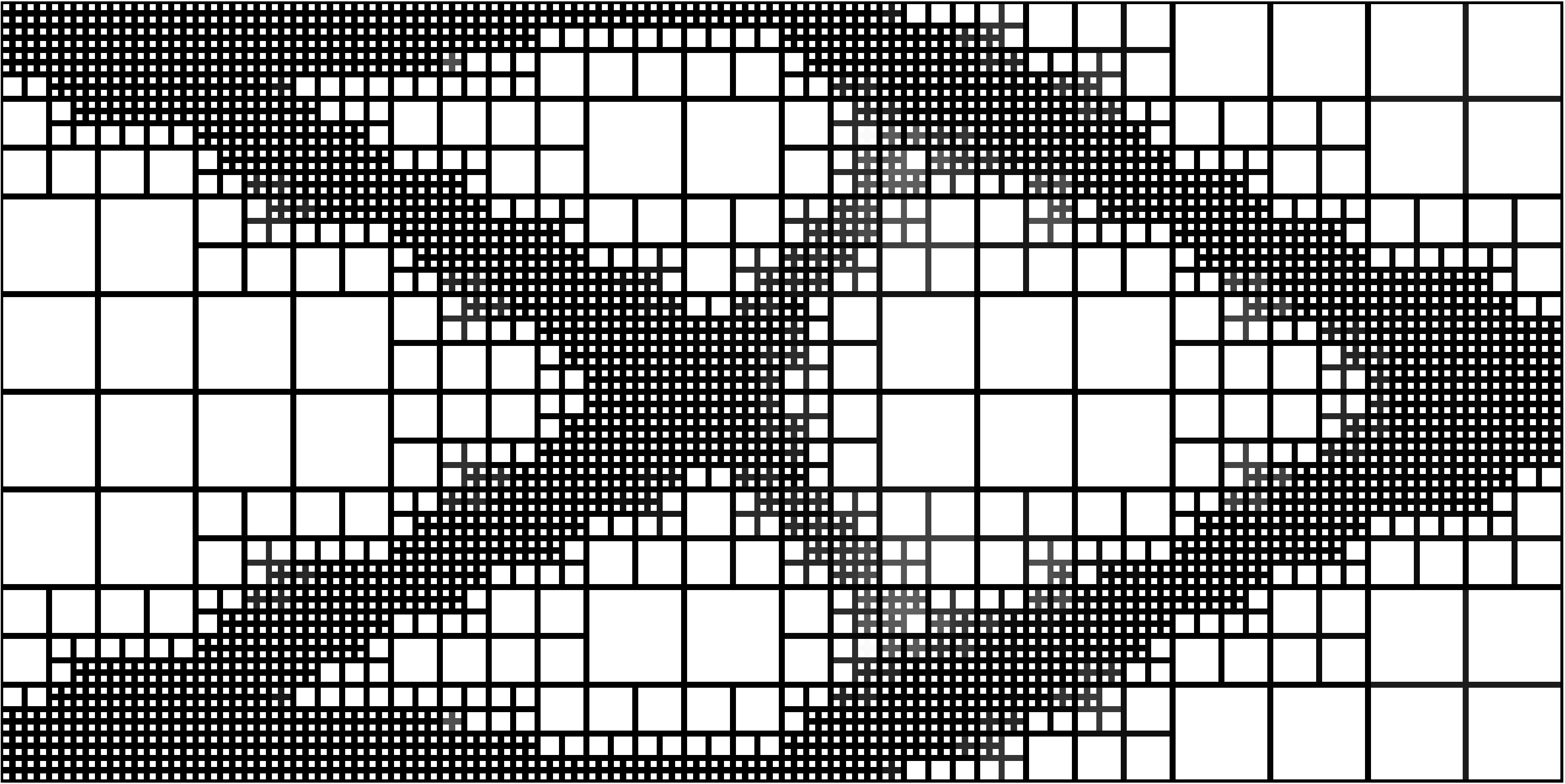}}
\caption{A sequence of the balanced quadtree during the optimization process. The number of iterations (It.), the compliance value ($c$), and the sharpness ($s$) are reported. 
}\label{fig:sequence}
\vspace{4mm}
\centering
\includegraphics[width=0.32\linewidth]{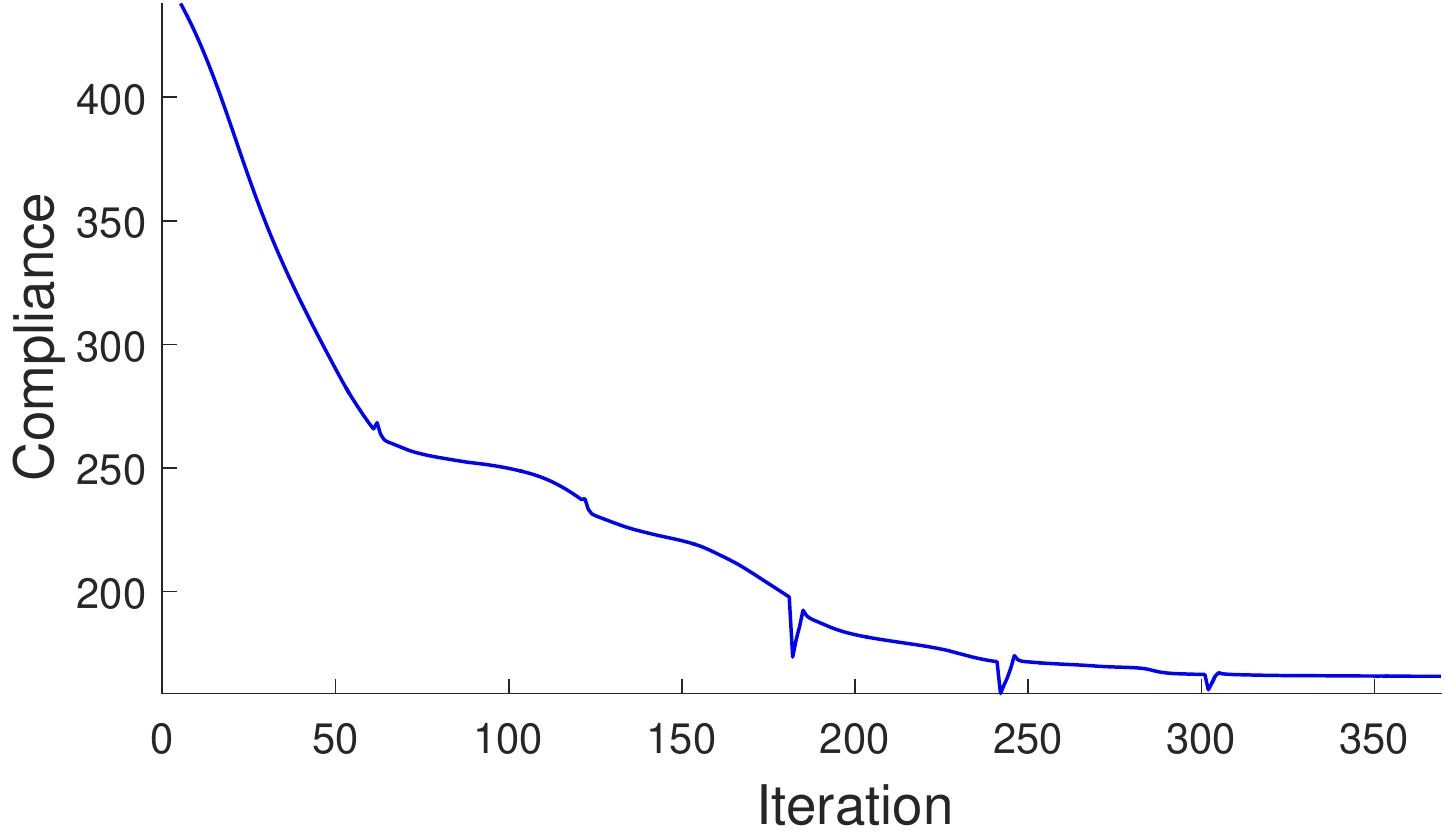} 
\includegraphics[width=0.32\linewidth]{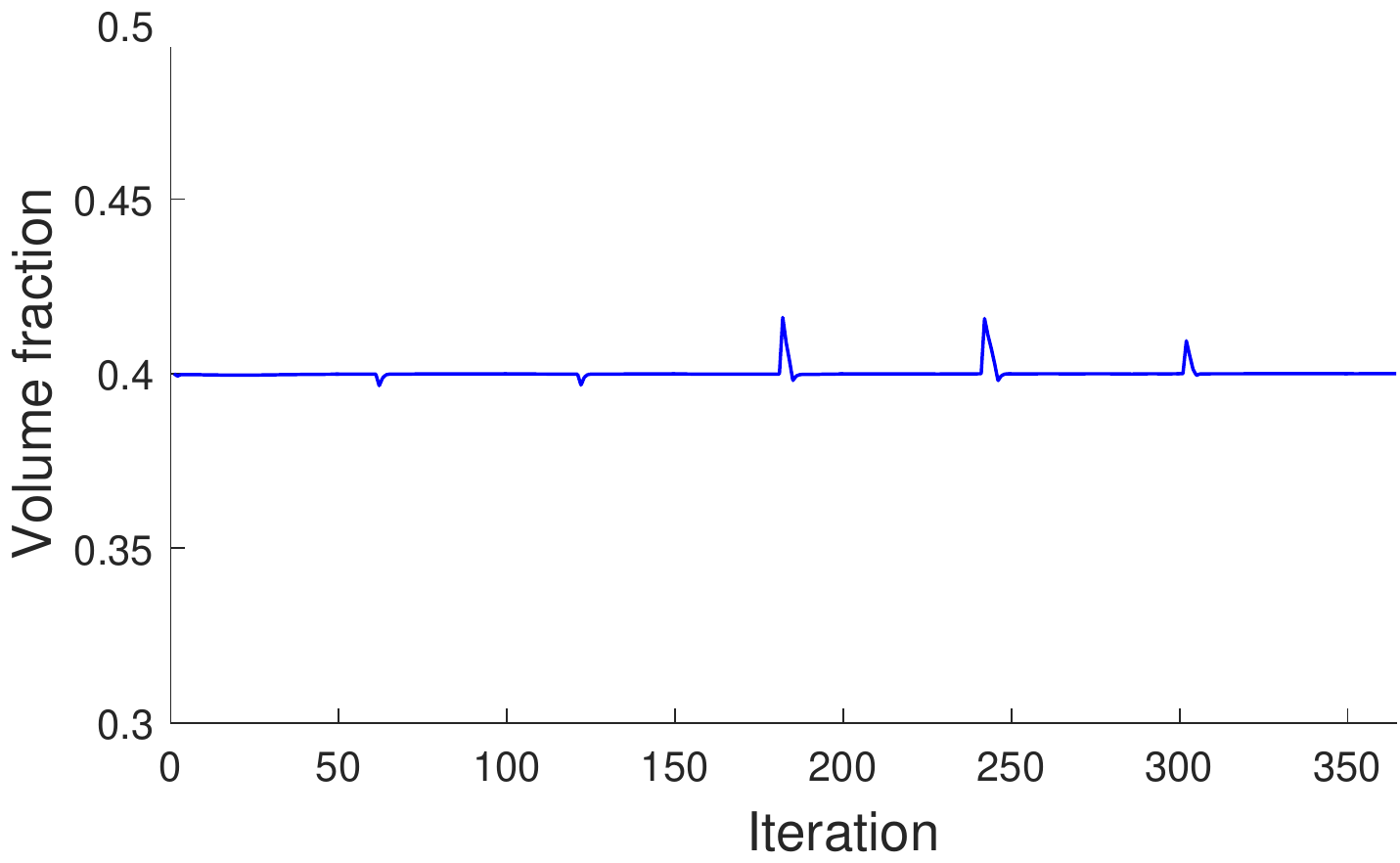}  
\includegraphics[width=0.32\linewidth]{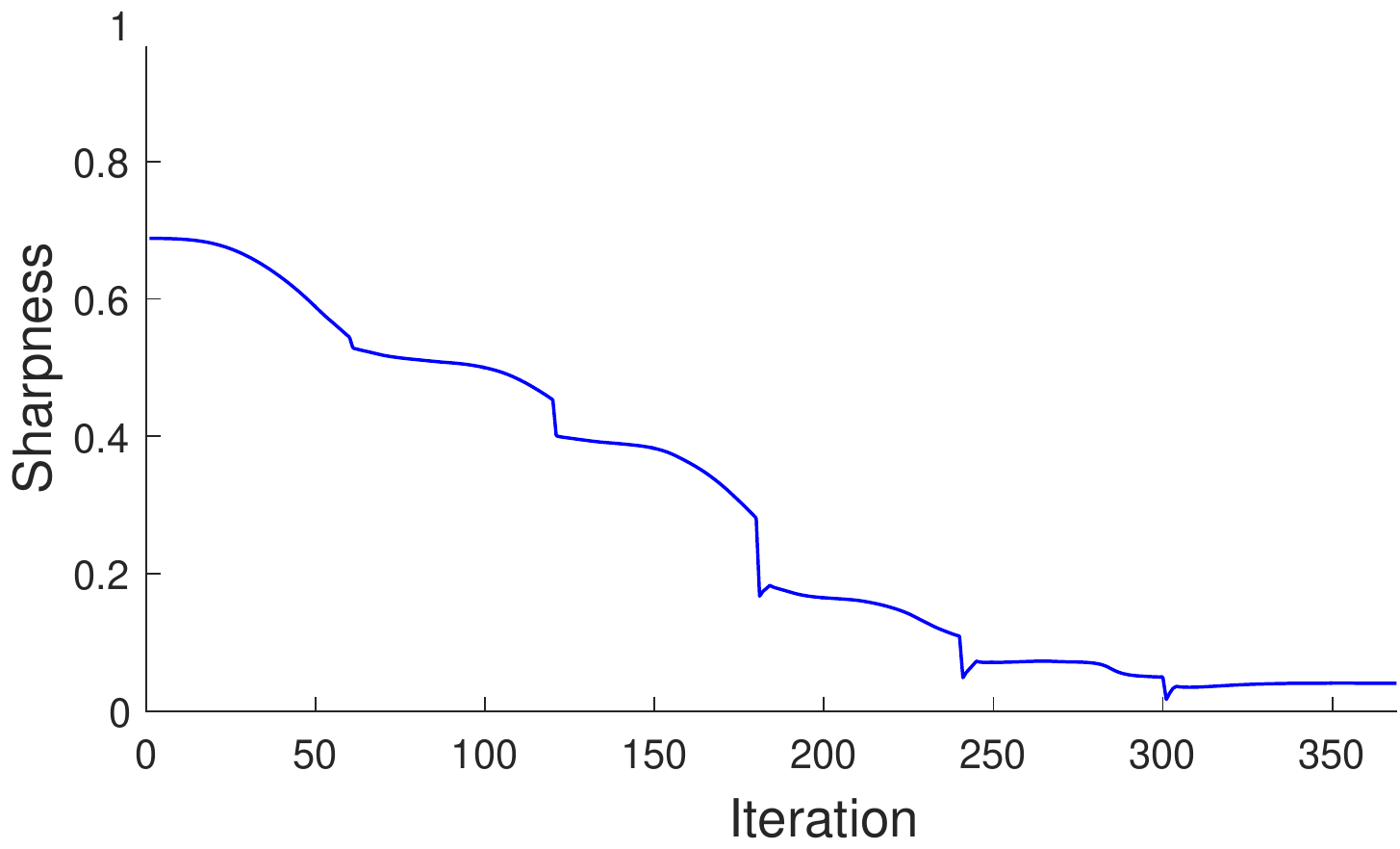}  
\caption{The convergence plots for the cantilever beam design (Fig.~\ref{fig:sequence}). From left to right: Convergence plot of compliance, volume fraction, and sharpness. The jumps in curves are due to the $\beta$-continuation in the Heaviside projection.
}\label{fig:Convergence}
\end{figure*}

The method has been implemented based on the Matlab code provided by Andreassen et al.~\cite{Andreassen10}. The optimization problem is solved using the method of moving asymptotes (MMA)~\cite{Svanberg87} \revise{}{and optimality criteria (OC)}. \revise{}{In our tests MMA yields a compliance value that is about $2\%$ smaller than the optimized value from OC. MMA results are reported in this paper.} Sensitivity analysis which is used in numerical optimization is given in Appendix.

In the following we report numerical analysis on simple design domains in Section~\ref{subsec:Analysis}, compare results from continuous and discrete optimization in Section~\ref{subsec:Comparison}, and discuss physical tests in Section~\ref{subsec:3Dprint}. The same parameters are used in all tests: The SIMP penalization $p=3$, the $p$-norm $p_n=-16$, the threshold for projection $\eta=0.5$. The sharpness of projection $\beta$ starts with $1$ and is increased up to 32, by doubling its value every $60$ iterations. The maximum number of iterations is $400$.

\subsection{Numerical Analysis}\label{subsec:Analysis}

\subsubsection{Cantilever Beam}\label{subsubsec:Cantilever}

The first example is a cantilever beam that has appeared multiple times during the exposition of the method. The boundary condition is illustrated in Fig.~\ref{fig:beamNoFilter}. The rectangular design domain is initialized with a coarse grid of $8 \times 4$. Each grid cell contains $64 \times 64$ square elements, leading to a finite element resolution of $512 \times 256$ for the entire domain. The maximum refinement level is $4$.

Figure~\ref{fig:sequence} shows a sequence of intermediate structures during the optimization of a balanced quadtree structure. As the optimization progresses, the quadtree structure emerges from the grey density distribution. To measure how close the continuous density field is to a binary field, the sharpness factor~\cite{Sigmund07} is defined as
\begin{equation}
s=\frac{4}{n}\sum_e(\rho_e (1-\rho_e)),
\end{equation}
where $n$ is the number of finite elements. The sharpness factor becomes $0.0$ when the density values converge to a strict 0-1 solution, and it becomes $1.0$ if all elements take a value of $0.5$. The optimized quadtree (Fig.~\ref{fig:sequence} bottom right) has a sharpness factor of $0.04$, which corresponds the field where density values on average take $0.01$ or $0.99$.

The convergence plots are shown in Fig.~\ref{fig:Convergence}. The plot on the left shows that the compliance value decreases almost monotonically, except for a few jumps corresponding to the increase of $\beta$ in the Heaviside projection. The middle shows the volume fraction during the process. The design variables are initialized homogeneously with a value of $0.49$ such that the volume fraction of $0.4$ has been satisfied from the beginning. The plot on the right shows the convergence of sharpness. The optimization is terminated at the iteration of $379$ as the maximum change in the density field becomes smaller than $1\mathrm{e}{-4}$.

\begin{figure}[!htb]
\centering
\def\svgwidth{0.98\linewidth}
\footnotesize	
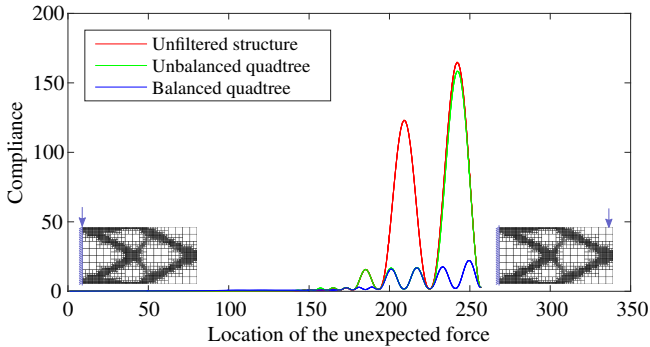
\caption{The three optimized structures (cf. Fig.~\ref{fig:beamNoFilter}, Fig.~\ref{fig:quadtreeBeamUnbalanced} and Fig.~\ref{fig:balancedBeam}) are tested under a small unexpected force, the location of which moves along the top boundary edge from left to right. The balanced quadtree is most robust as the corresponding curve (blue) shows the smallest oscillation.
}\label{fig:robust}
\end{figure}

One motivation to have quadtree structures is their potential robustness in handling unexpected perturbation forces which are not modelled in numerical optimization. Fig.~\ref{fig:robust} shows a comparison on the three optimized structures: the unfiltered structure (Fig.~\ref{fig:beamNoFilter}), the unbalanced quadtree (Fig.~\ref{fig:quadtreeBeamUnbalanced}), and the balanced version (Fig.~\ref{fig:balancedBeam}). In this test, we simulate a small, unexpected force. Its location moves along the top boundary edge, while its magnitude is one tenth of the force which was used in the deterministic structural optimization. The fixation on the left boundary edge is unaltered in this test, but it shall be noted that the fixation is also subject to uncertainty depending on the use scenario. The first half of the three curves is very stable. As the force moves to the right side where the structure becomes sparse, we observe oscillations in compliance curves corresponding to the periodicity of the underlying quadtree cells. The blue curve has a oscillation significantly smaller than the other two, meaning that the balanced quadtree is most robust. 

\subsubsection{MBB Beam}\label{subsubsec:MBB}

\begin{figure}[!htb]
\centering
\small
\def\svgwidth{0.98\linewidth}
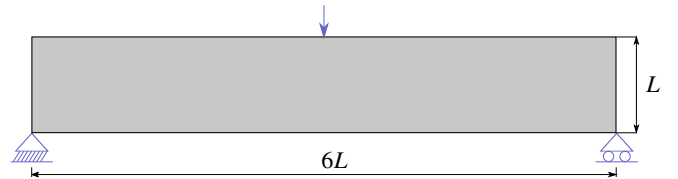
\caption{Illustration of the MBB beam design problem.
}\label{fig:MBBBC}
\end{figure}

The second example is the MBB (Messerschmitt-B{\"o}lkow-Blohm) design problem shown in Fig.~\ref{fig:MBBBC}. Due to symmetry only the right half of the design domain is optimized. The prescribed material volume is $30\%$ of the design domain. The half domain is initialized with a coarse grid of $12\times4$, where each grid cell is composed of $64\times64$ elements, leading to a total of $768\times256$. 

\begin{figure*}[!htb]
\centering
\subfloat[Unbalanced quadtree, $\bar{k}=3$, $c=1083.5$, $s=0.002$]{\includegraphics[width=0.4\textwidth]{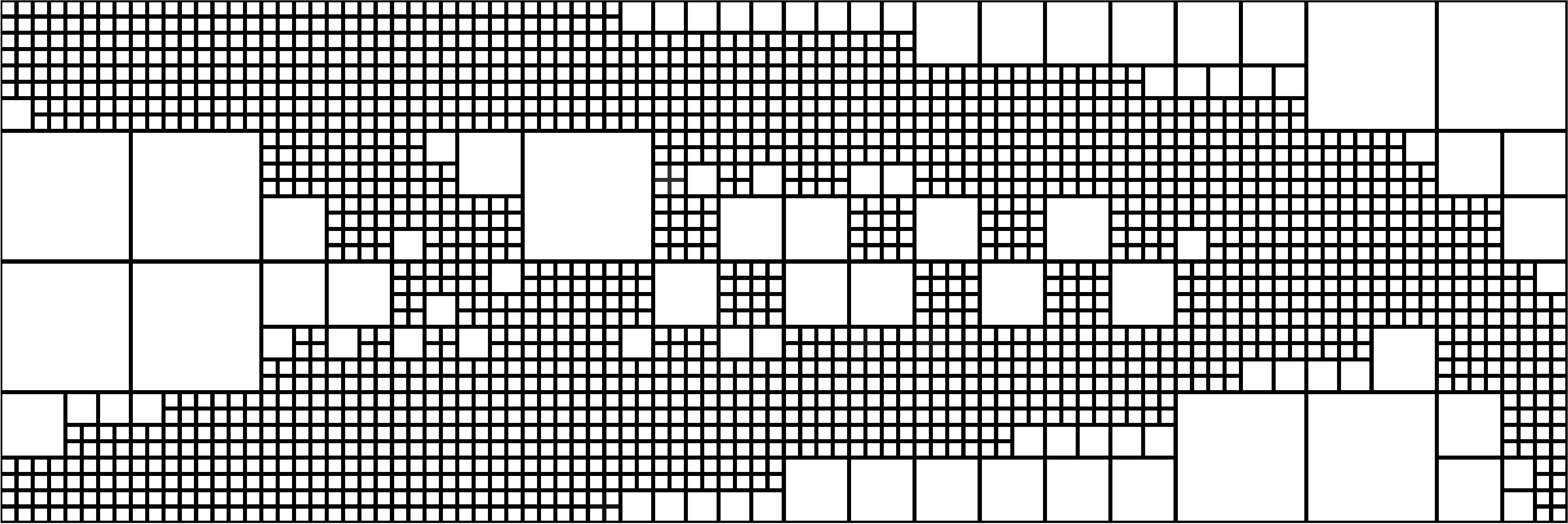}}
\qquad
\subfloat[Balanced quadtree, $\bar{k}=3$, $c=1127.9$, $s=0.008$]{\includegraphics[width=0.4\textwidth]{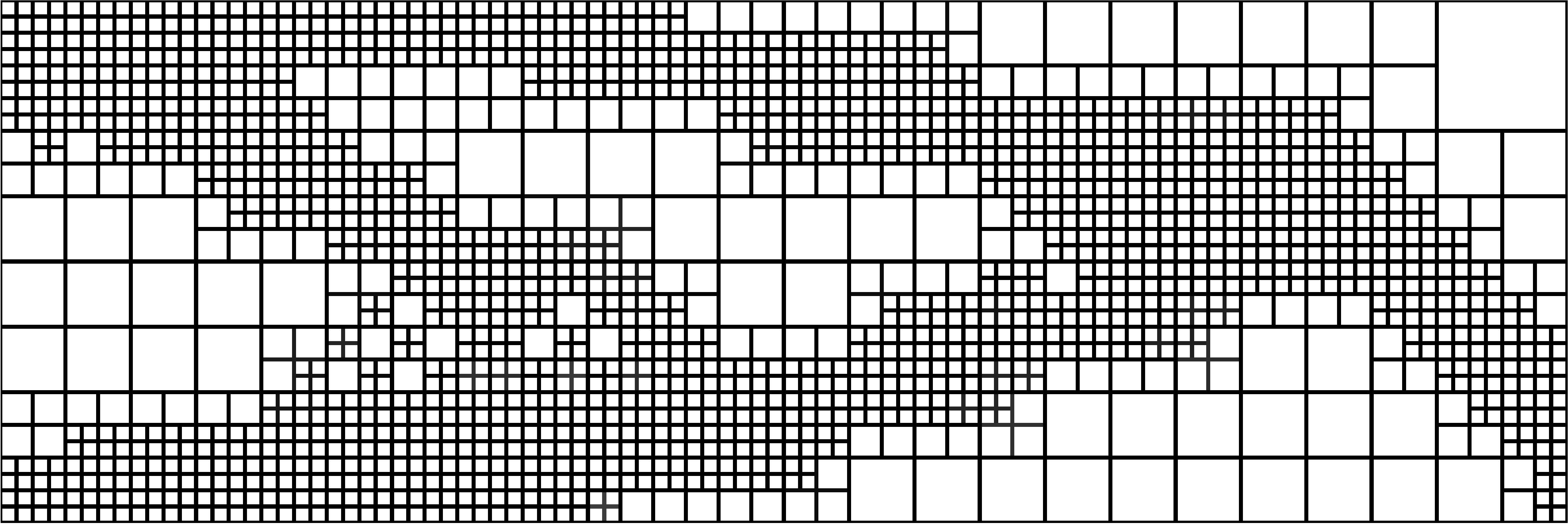}}
\\
\subfloat[Unbalanced quadtree, $\bar{k}=4$, $c=571.6$, $s=0.010$]{\includegraphics[width=0.4\textwidth]{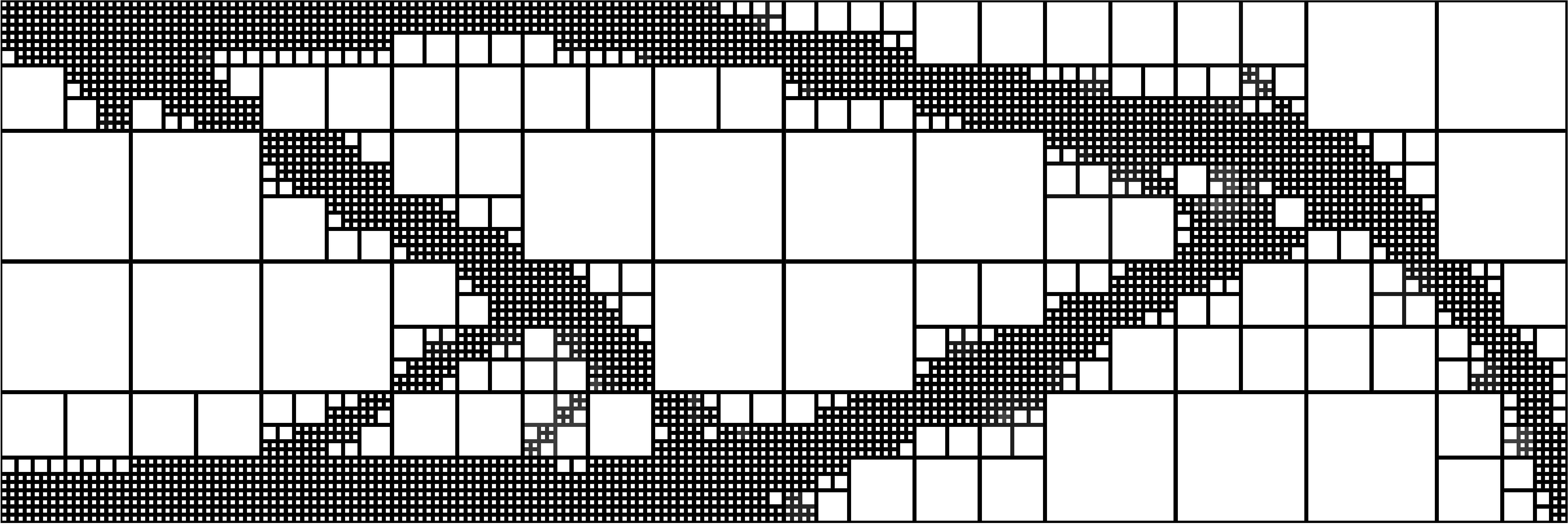}} 
\qquad
\subfloat[Balanced quadtree, $\bar{k}=4$, $c=711.5$, $s=0.028$]{\includegraphics[width=0.4\textwidth]{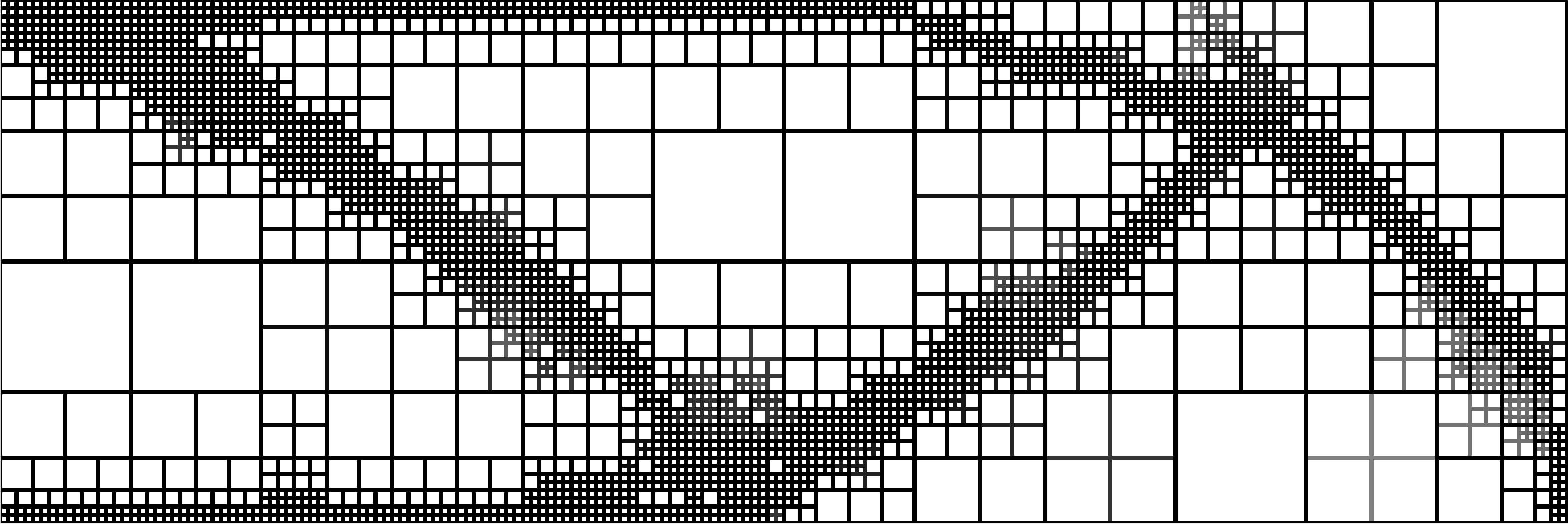}}
\\
\subfloat[Unbalanced quadtree, $\bar{k}=5$, $c=449.0$, $s=0.030$]{\includegraphics[width=0.4\textwidth]{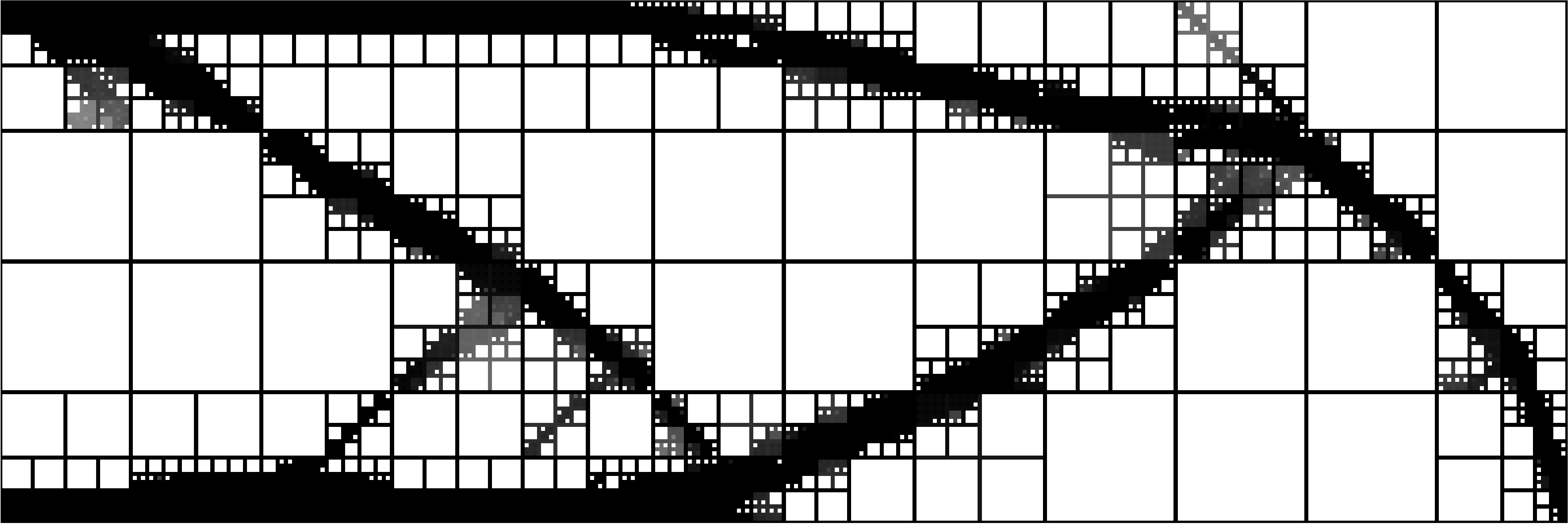}} 
\qquad
\subfloat[Balanced quadtree, $\bar{k}=5$, $c=586.1$, $s=0.070$]{\includegraphics[width=0.4\textwidth]{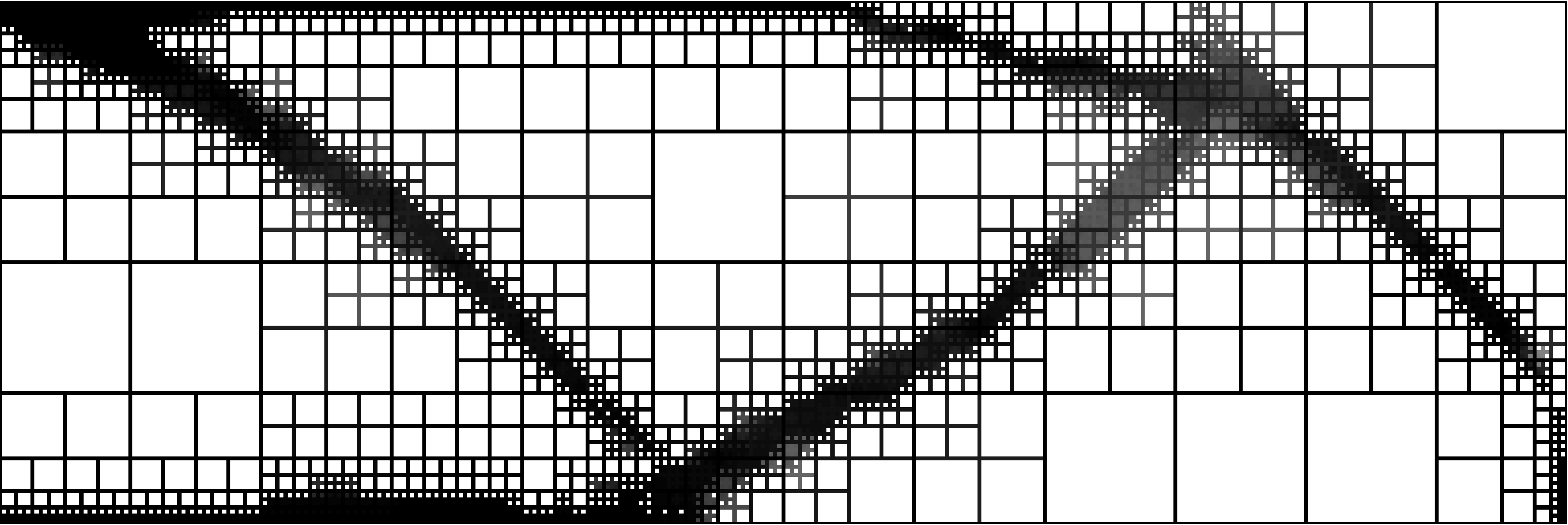}}
\caption{A comparison of unbalanced (left) and balanced quadtree (right) using a refinement level of, from top to bottom, three, four, and five.
}%
\label{fig:MBB}%
\vspace{4mm}
\centering
\includegraphics[width=0.32\linewidth]{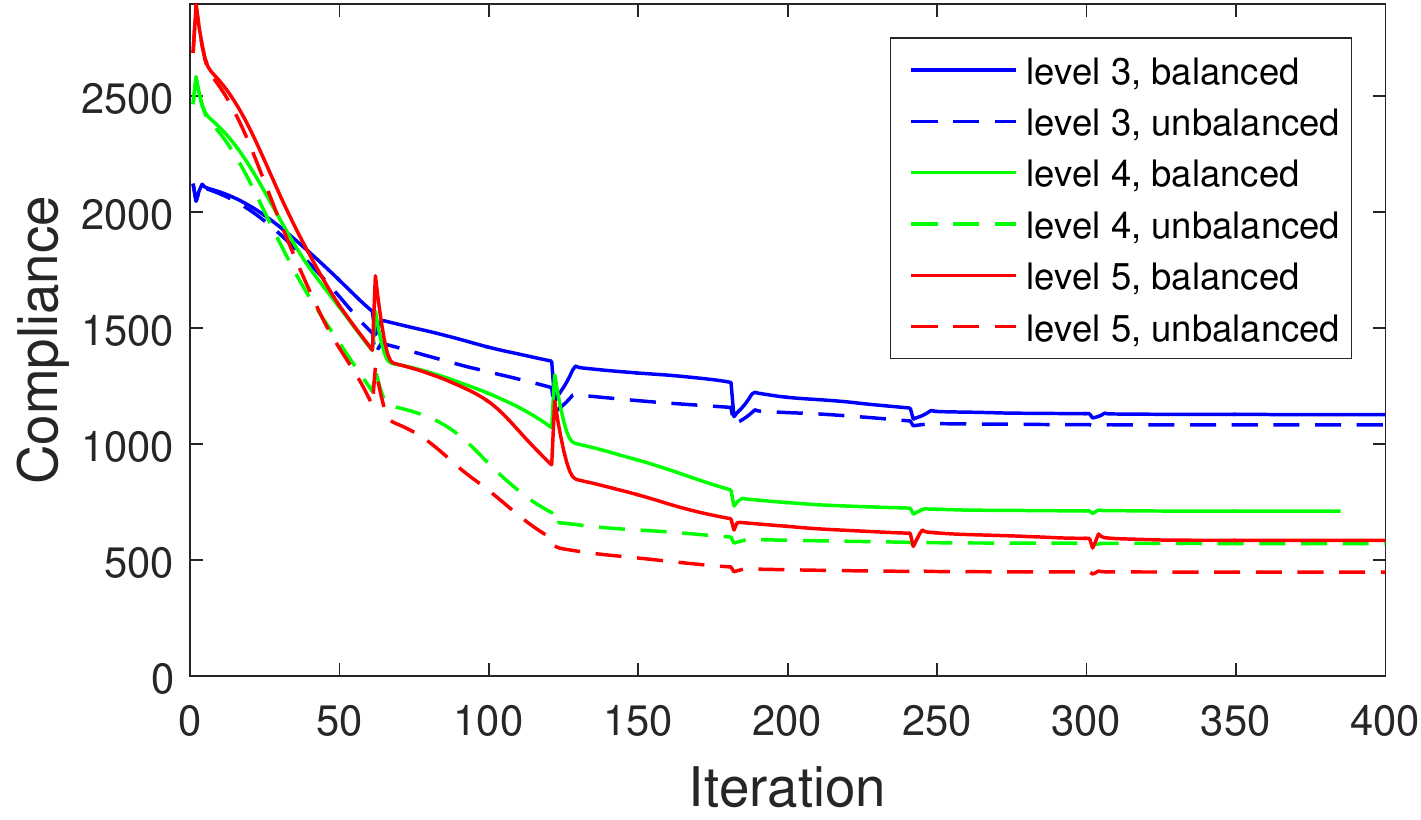} 
\includegraphics[width=0.32\linewidth]{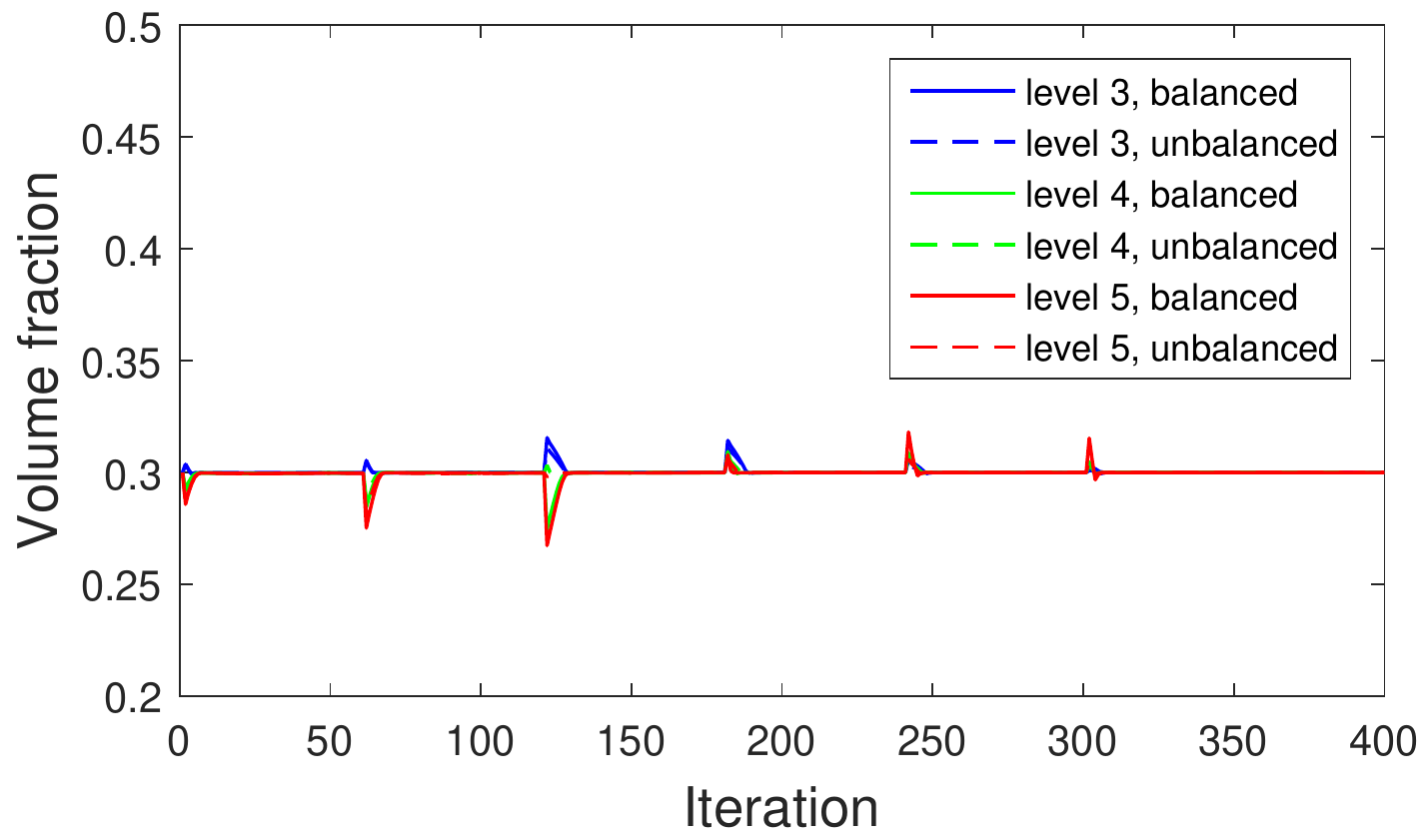}  
\includegraphics[width=0.32\linewidth]{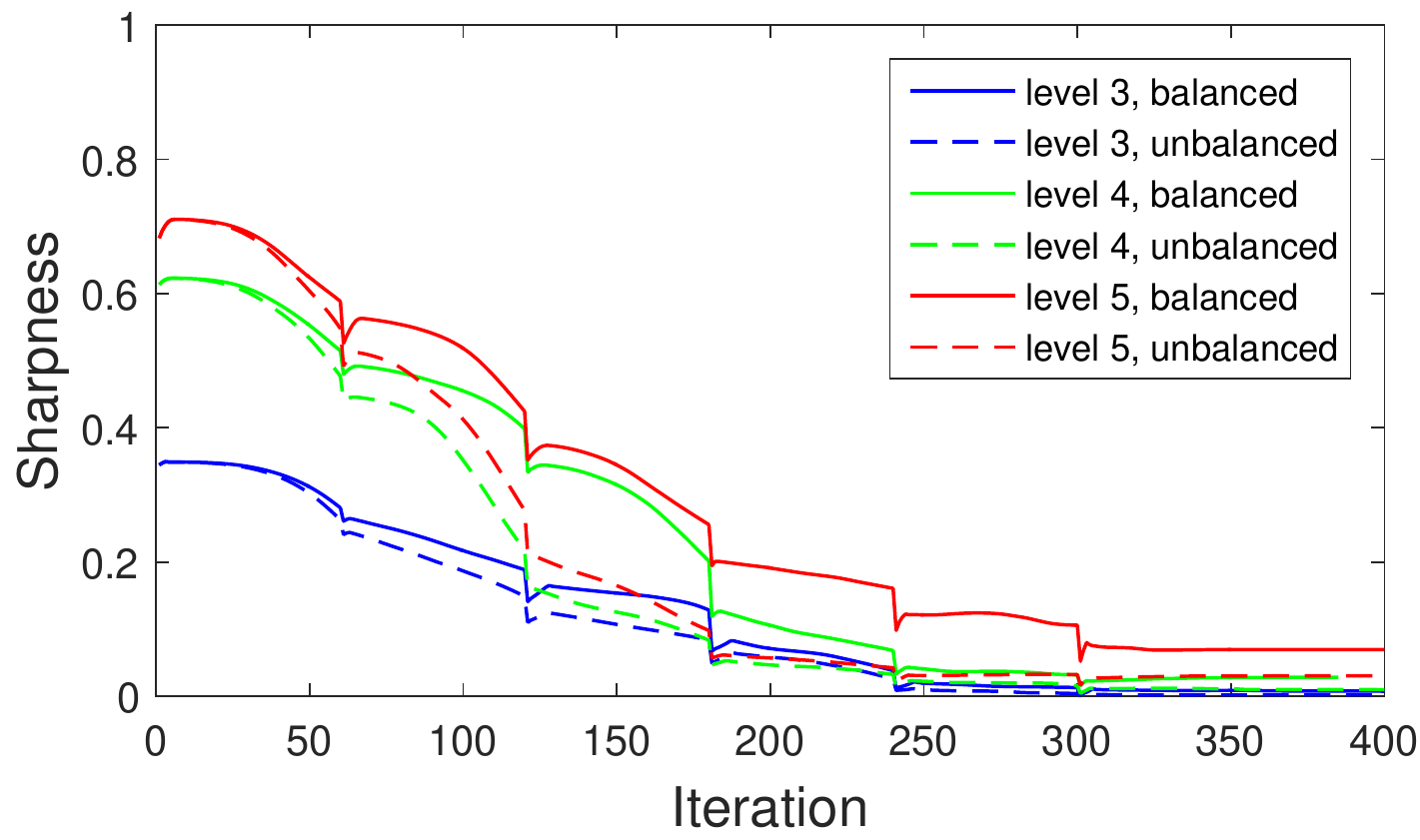}  
\caption{From left to right: The convergence of compliance, volume fraction, and sharpness for the MBB designs shown in Fig.~\ref{fig:MBB}. 
}\label{fig:MBBConvergence}
\end{figure*}

Figure~\ref{fig:MBB} compares the quadtree structures obtained using different refinement levels. The three rows from top to bottom show optimized quadtree using a maximum refinement level of three, four, and five, respectively. In each row, the left corresponds to the unbalanced refinement, while the right is obtained in the setting of balanced refinement. Comparing left and right, it can be found that the balanced quadtree (right) has a larger compliance and a larger sharpness value than in the unbalanced setting (left). This is due to the fact that the balanced refinement additionally involves the neighbours of the parent cells, i.e., more restrictions on refinement. Comparing the three rows, we observe smaller compliances when a larger refinement level is applied. This is because a larger refinement level increases the solution space. At a refinement level of five (bottom), small cavities from the previous refinement (middle) are filled with material. The optimizations with these different settings all converge well, as can be seen from the convergence plots in Fig.~\ref{fig:MBBConvergence}.

\begin{figure*}[!htb]
\centering
\subfloat[Boundary condition]{\def\svgwidth{0.3\linewidth}
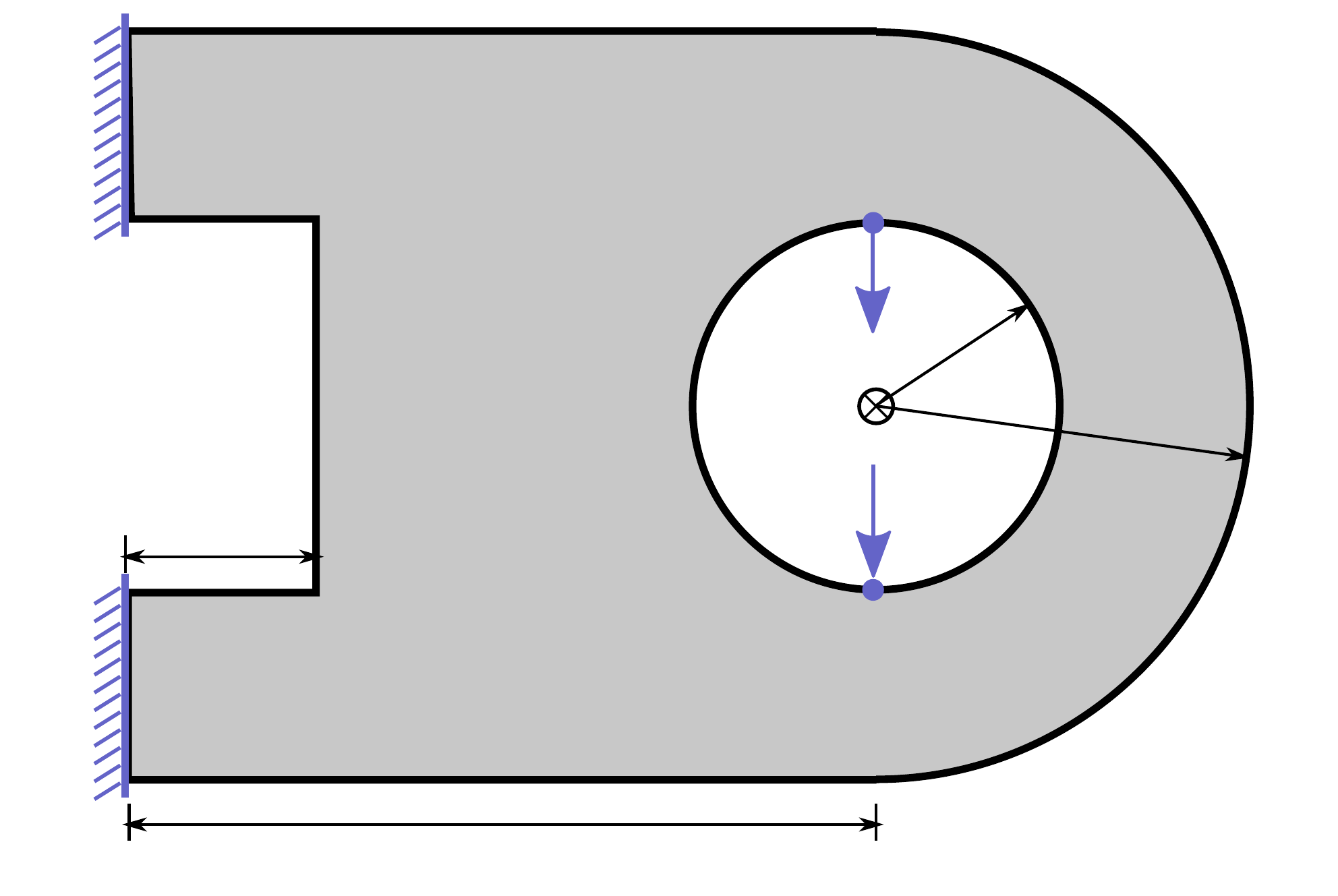}
\subfloat[Uniform pattern ($c=566.2$)]{\includegraphics[width=0.3\linewidth]{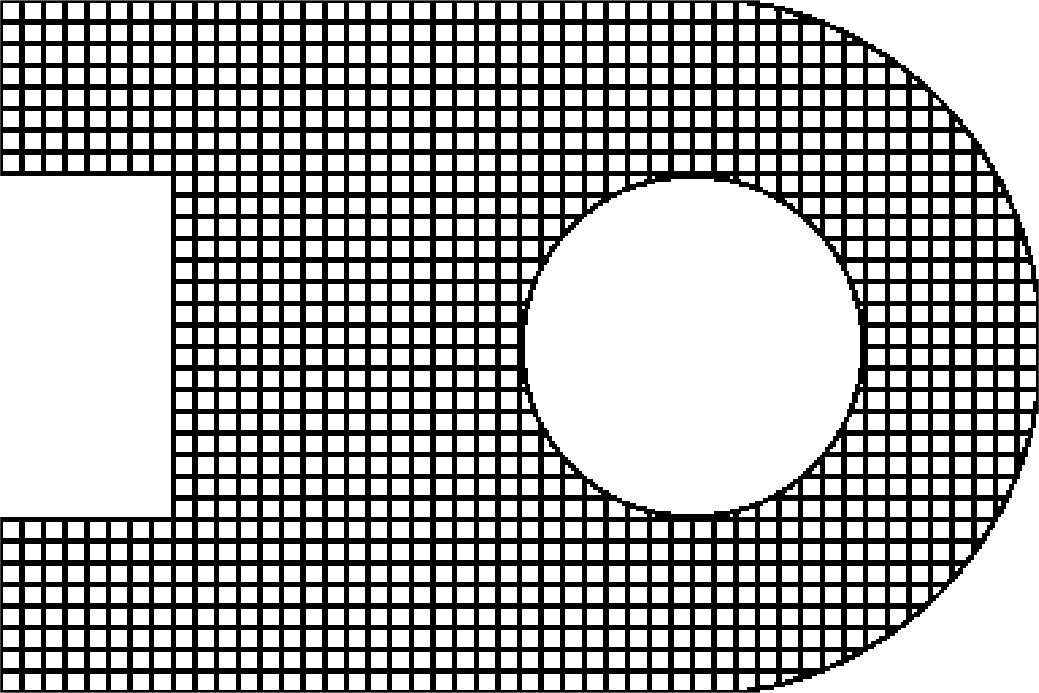}}
\qquad
\subfloat[Balanced quadtree ($c=143.4$)]{\includegraphics[width=0.3\linewidth]{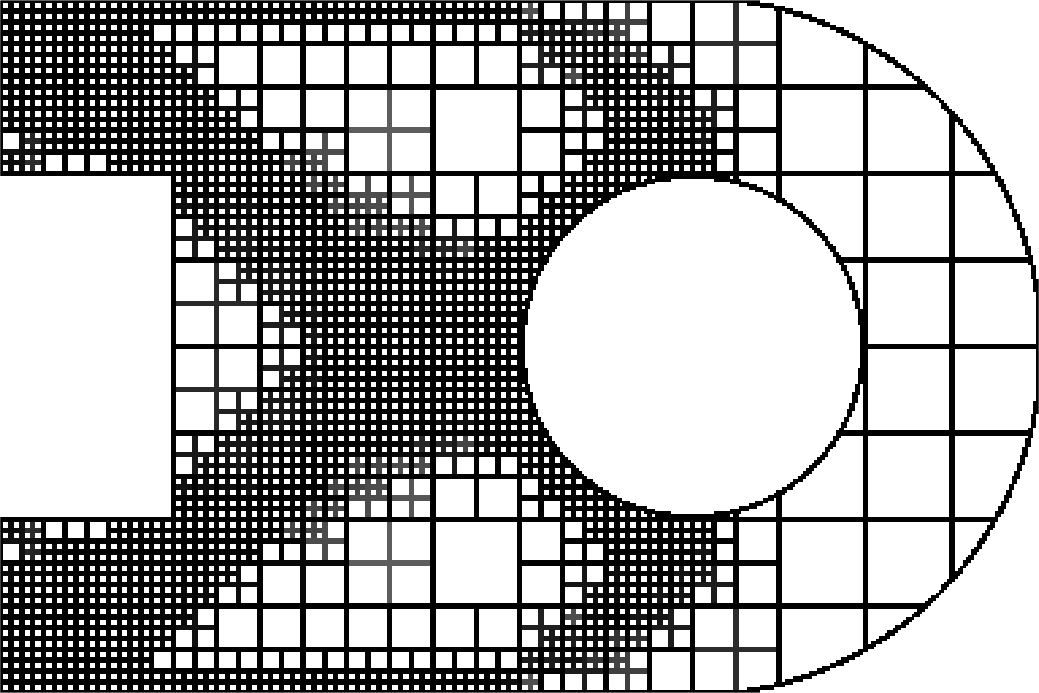}}
\caption{Quadtree optimization on a curved design domain. The adaptively refined quadtree is $3.9$ times stiffer than the uniform pattern. 
}\label{fig:bracket2}
\end{figure*}

\subsubsection{Bracket}\label{subsubsec:Bracket}
The adaptive structure refinement works as well on curved design domains. For a curved design domain, we simulate the smallest axis-aligned bounding box enclosing the curved domain. The voids outside the curved domain are assigned as passive void elements, while a boundary layer with a thickness of 2 elements are passive solid elements. Passive elements are excluded from the design update. This requires modifying the corresponding entries in the transformation matrix in Eq.~\eqref{eq:MappingFilter}. 


\begin{figure*}[t]
\centering
\small
\def\svgwidth{0.9\linewidth}
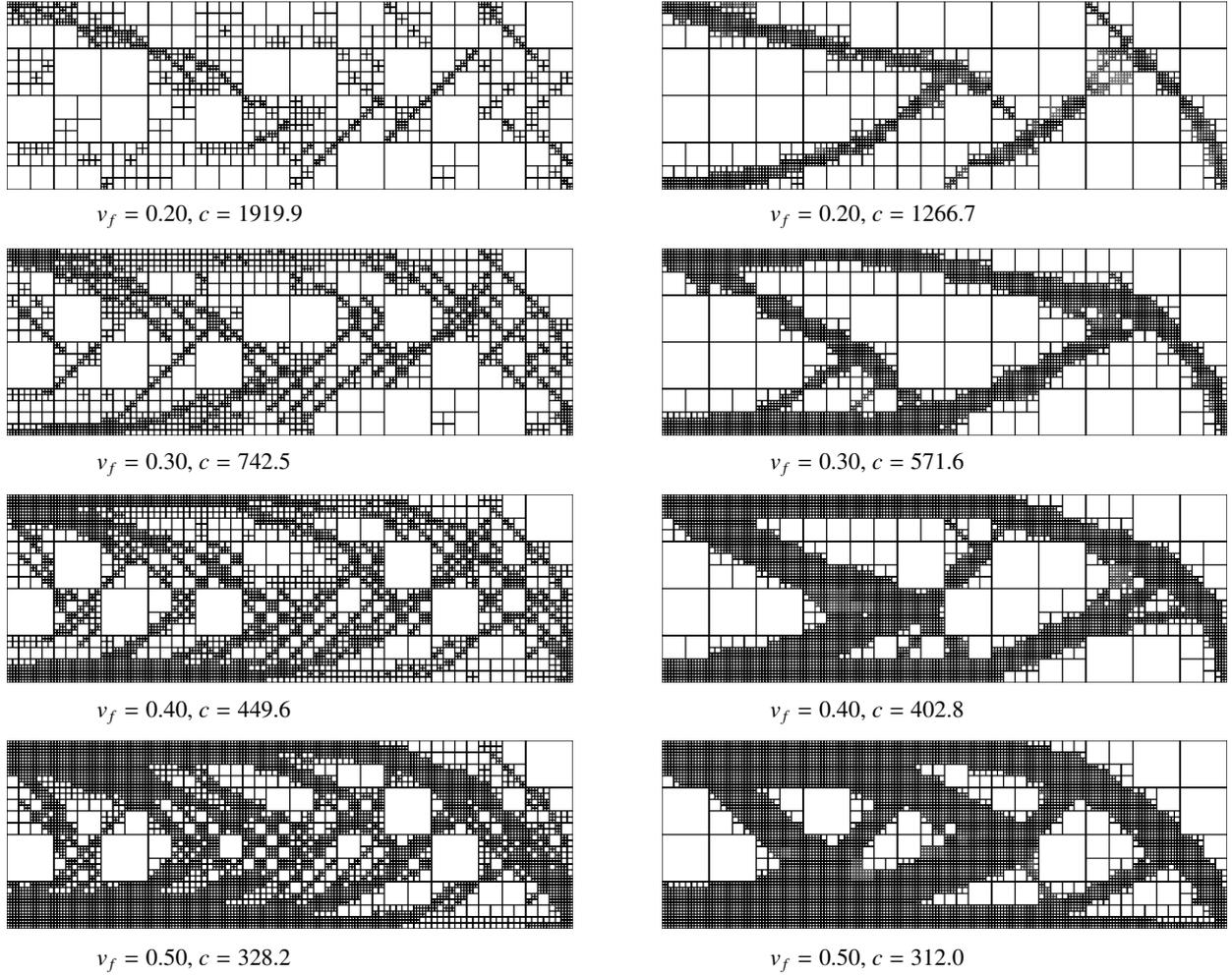
\caption{Quadtree structures obtained by a greedy approach in the discrete formulation (left) and optimized in the continuous formulation (right). From top to bottom the comparisons are made under an increasing volume fraction from $0.2$ to $0.5$.
}\label{fig:comparisonBinary}
\end{figure*}

Figure~\ref{fig:bracket2} (a) illustrates the dimension and boundary conditions of the bracket design problem. It contains curved boundary where a quadtree mesh is not conformal. The tight bounding box of the bracket is discretized using a coarse grid of $12\times 8$, where each grid cell has $32\times 32$ elements, i.e., a total of $384\times 256$ elements for the bounding box. First we uniformly refine the coarse grid by two level which is shown in Fig.~\ref{fig:bracket2}~(b). Its volume fraction is $33.02\%$. Using this volume fraction as a constraint we perform adaptive structure optimization, allowing three levels of refinement. The optimized quadtree is shown in Fig.~\ref{fig:bracket2}~(c). Numerical analysis suggests that the adaptive quadtree is about 4 times stiffer than the uniform pattern. 

\subsection{Comparison to Greedy Approach}\label{subsec:Comparison}

The quadtree structures optimized by the continuous formulation are compared with those obtained by the greedy approach~\cite{Wu16Rhombic}. The greedy approach starts refining from the initial coarse grid, and performs selective refinement at odd iterations and selective coarsening at even iterations. In each odd iteration, following finite element analysis, the leaf cells in the quadtree are sorted by their strain energy density. The first few cells in this sorted list are then successively refined, until the volume increase exceeds a small volume change, $0.4\%$ in our examples. Alternatively in each even iteration, following finite element analysis, blocks of $2^2$ leave cells, i.e., the cells on the second finest level, are sorted by the averaged strain energy density. The last few cells in this sorted list are coarsened, until the volume is decreased by $0.1\%$. The combination of refinement and coarsening is found more optimal than purely refinement, and thus serves a fair reference for comparison. The greedy approach is only applied in the setting of unbalanced quadtree.

\begin{figure*}[!htb]
\centering
\small
\def\svgwidth{0.9\linewidth}
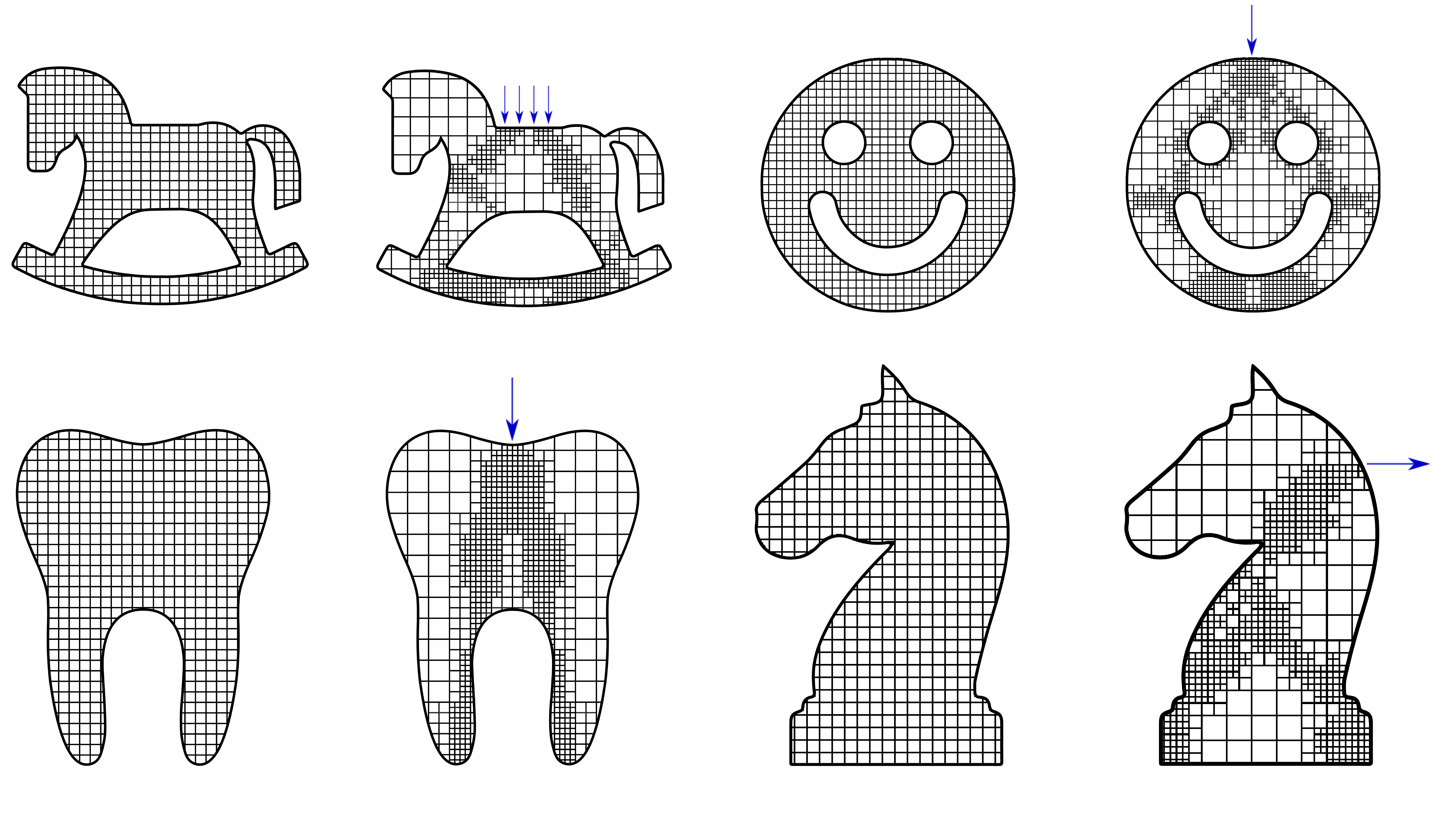\\
\vspace{5mm}
\includegraphics[width=0.9\linewidth]{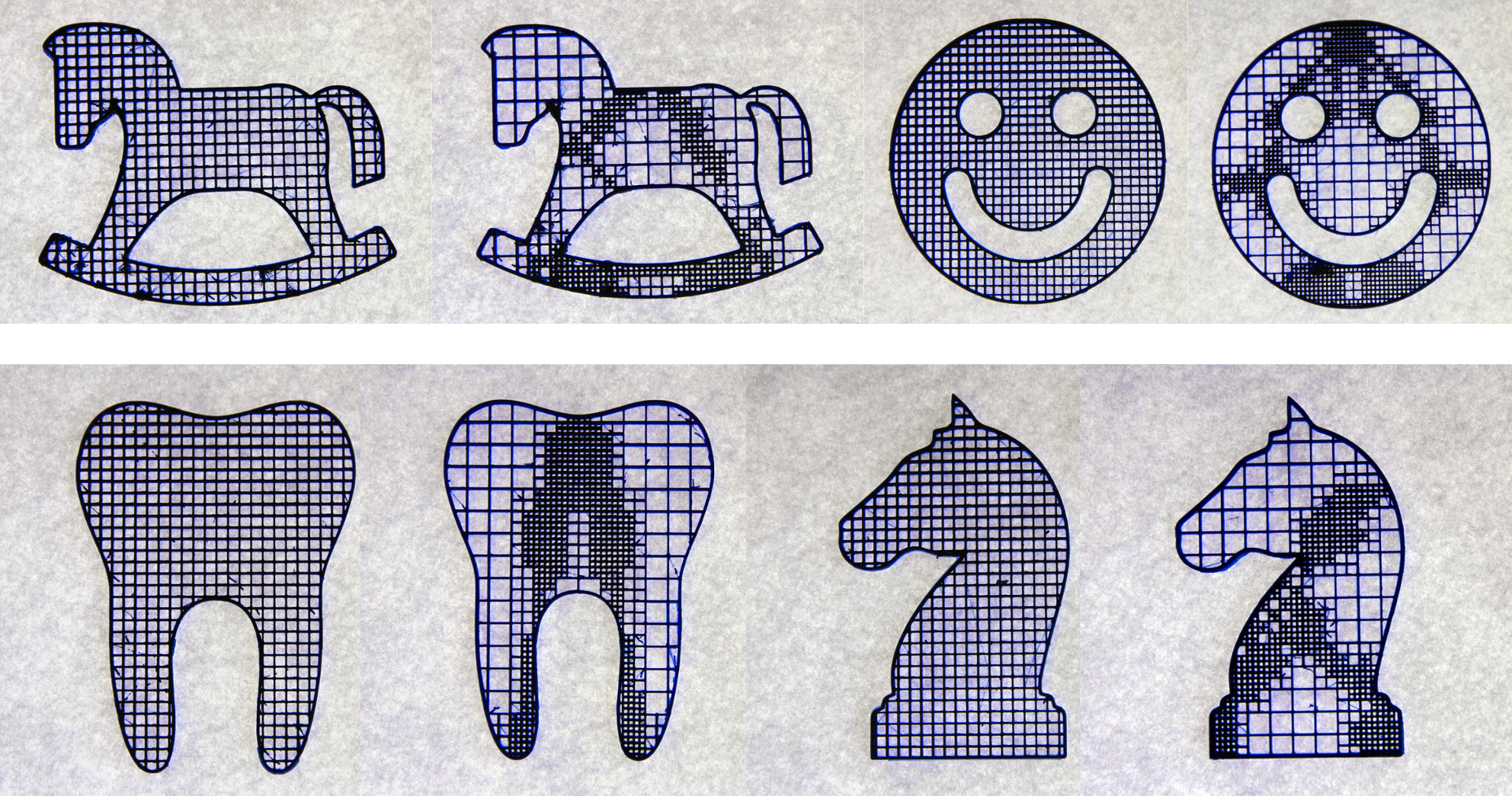}  
\caption{Digital structures (top two rows) and 3D printed models (bottom two rows).
}\label{fig:3dprint}
\end{figure*}

Figure~\ref{fig:comparisonBinary} shows the optimized MBB structures by the greedy approach (left) and by the continuous formulation (right). The structures are compared under four different volume fractions. The structures in the right column show some distinct dense and sparse regions, to some extent similar to classic topology optimization, while structures in the left column have regions of alternating refined and non-refined cells, and the alternation visually appears a bit random. The quadtree structures optimized by the continuous formulation consistently have a smaller compliance value. The difference in compliance values is especially significant when the volume fraction is small (cf. the top row). \revise{}{The greedy approach takes  $89$, $160$, $229$, and $297$ iterations to achieve a volume fraction of $0.2$, $0.3$, $0.4$, and $0.5$, respectively, while in the continuous optimization a maximum of $400$ iterations for each volume fraction is prescribed.}

\subsection{Physical Test}\label{subsec:3Dprint}
We further perform physical tests to compare the stiffness of uniform patterns and optimized quadtree structures. To this end, we have 3D printed the structures with an Ultimaker~3, which uses the Fused Deposition Modelling (FDM) process. The printing material is PLA. The nozzle size is $0.4\,mm$.  Fig.~\ref{fig:3dprint} shows the 3D printed models. For each pair, the volume taken by the uniform pattern is prescribed as the volume limit for quadtree optimization. The 2D models are scaled to fit within a domain of $100\,mm \times 100\,mm$. They are then extruded to 3D with a thickness of $5\,mm$. 

\begin{figure}[t]
\centering
\includegraphics[width=0.98\linewidth]{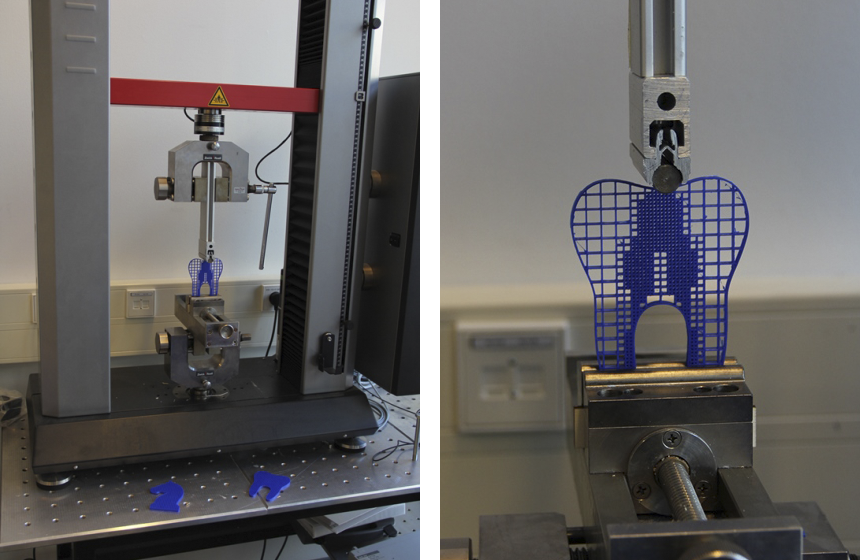}
\caption{Physical test on the fabricated model.
}\label{fig:testsetup}
\end{figure}

\begin{figure}[t]
\centering
\includegraphics[width=0.8\linewidth]{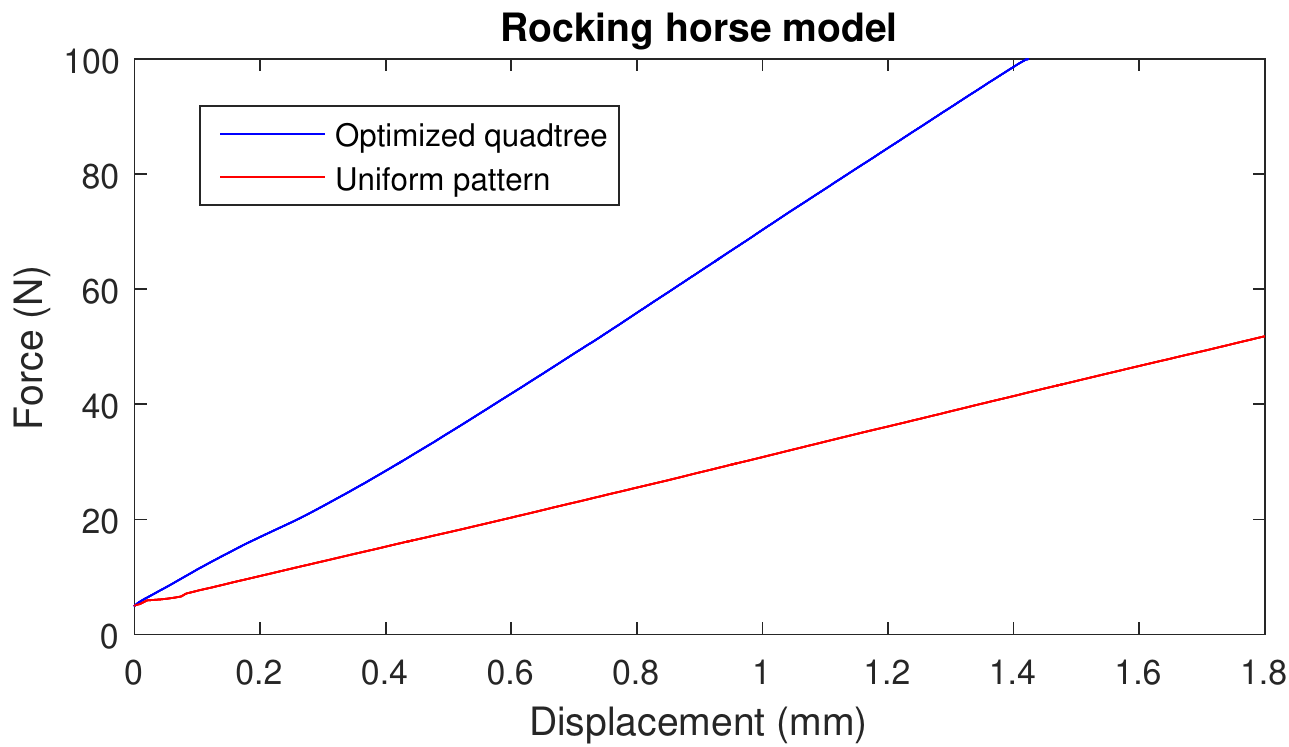}\\
\vspace{4mm}
\includegraphics[width=0.8\linewidth]{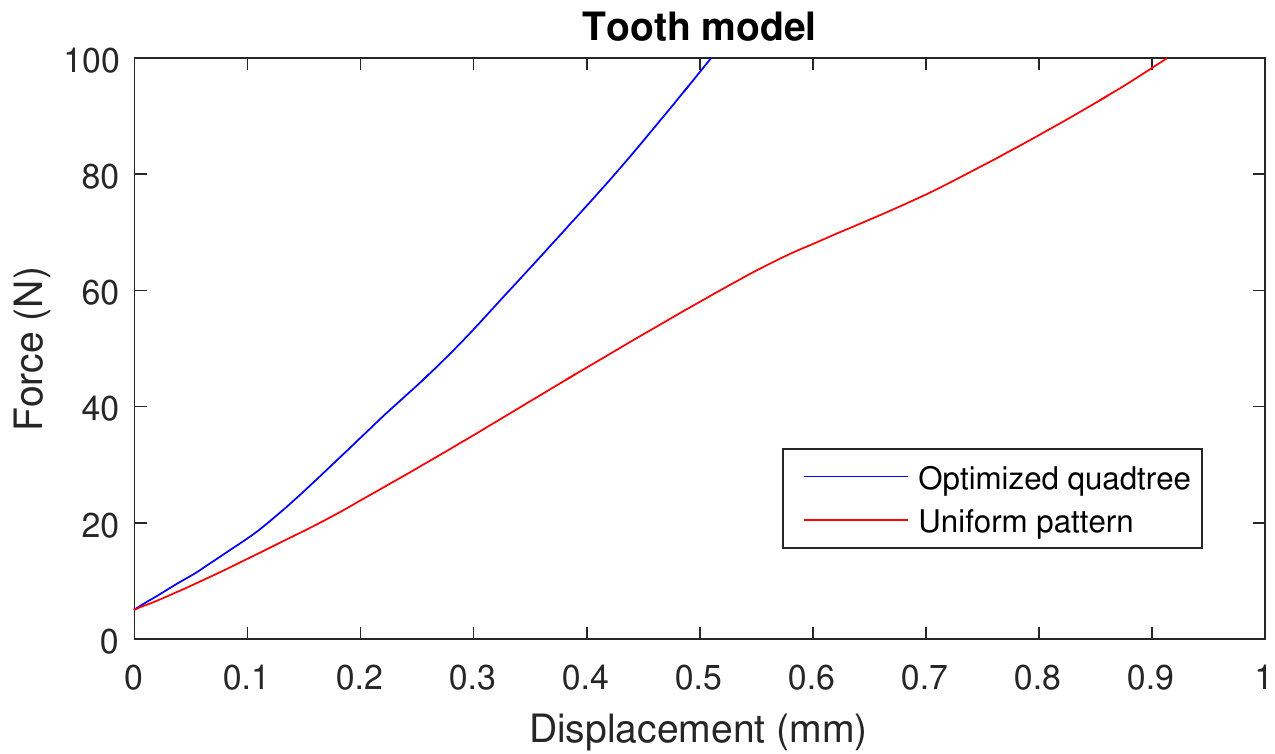}
\caption{The measured force-displacement curves on the rocking horse model (top) and the tooth model (bottom). 
}\label{fig:tensileResults}
\end{figure}

The experimental setup is shown in Fig.~\ref{fig:testsetup}. The models are supported from beneath, while the clamp is used to keep their position. During the physical test the models are pre-loaded with a $5$ Newton force. The force then increases at a rate of $10\,N/s$ until it reaches $100$ Newton. Fig.~\ref{fig:tensileResults} shows the force-displacement curves for the rocking hose and tooth models. In both tests, the curve corresponding to the optimized quadtree is much steeper, confirming that the optimized quadtree has a much higher stiffness. 

\section{Conclusion}\label{sec:Conclusion}

We have presented a continuous optimization method to design adaptively refined structures. The structures are refined within a closed-walled design domain. The adaptive structure offers highly optimized stiffness for prescribed mechanical loads, and outperforms the uniform patterns.
Meanwhile the adaptive structure spreads over the interior of the mechanical part with controllable and smoothly varying void sizes, making it stable for uncertain and small perturbation forces. Continuous optimization achieves solutions that are more optimal compared to the greedy approach for solving the discrete formulation. \revise{}{Additional tests confirm our method works as well for multiple load cases and are left out for space reasons.}

Our method for continuous quadtree optimization opens up multiple directions. First, extending from quadtree to octree is straightforward. In 3D the structural elements can be either frames or walls, mapped from edges or faces of the octree. \revise{}{3D rhombic structures that are refined using a greedy approach have been reported in~\cite{Wu16Rhombic}.} Second, it is very interesting to design hierarchical structures in other forms, e.g., triangles/hexagons, or even fractals. Third, to alleviate the intensive computation in finite element analysis (especially for 3D), truss/beam elements and numerical homogenization will be useful.

\section*{Acknowledgements}
We thank Rob Scharff for preparing the physical tests, and Krister Svanberg for the Matlab MMA code.

\section*{References}

\bibliographystyle{elsarticle-num}
\bibliography{top.bib}

\begin{thebibliography}{10}
\expandafter\ifx\csname url\endcsname\relax
  \def\url#1{\texttt{#1}}\fi
\expandafter\ifx\csname urlprefix\endcsname\relax\def\urlprefix{URL }\fi
\expandafter\ifx\csname href\endcsname\relax
  \def\href#1#2{#2} \def\path#1{#1}\fi

\bibitem{Rajeev1992}
S.~Rajeev, C.~Krishnamoorthy, Discrete optimization of structures using genetic
  algorithms, Journal of structural engineering 118~(5) (1992) 1233--1250.

\bibitem{Wu16Rhombic}
J.~Wu, C.~C. Wang, X.~Zhang, R.~Westermann, Self-supporting rhombic infill
  structures for additive manufacturing, Computer-Aided Design 80 (2016) 32 --
  42.
\newblock \href {http://dx.doi.org/10.1016/j.cad.2016.07.006}
  {\path{doi:10.1016/j.cad.2016.07.006}}.

\bibitem{Sigmund13}
O.~Sigmund, K.~Maute, Topology optimization approaches, Struct. Multidiscip.
  Optim. 48~(6) (2013) 1031--1055.
\newblock \href {http://dx.doi.org/10.1007/s00158-013-0978-6}
  {\path{doi:10.1007/s00158-013-0978-6}}.

\bibitem{Deaton14}
J.~D. Deaton, R.~V. Grandhi, A survey of structural and multidisciplinary
  continuum topology optimization: post 2000, Struct. Multidiscip. Optim.
  49~(1) (2014) 1--38.
\newblock \href {http://dx.doi.org/10.1007/s00158-013-0956-z}
  {\path{doi:10.1007/s00158-013-0956-z}}.

\bibitem{Berger84}
M.~J. Berger, J.~Oliger, Adaptive mesh refinement for hyperbolic partial
  differential equations, Journal of Computational Physics 53~(3) (1984)
  484--512.

\bibitem{Maute95}
K.~Maute, E.~Ramm, Adaptive topology optimization, Structural optimization
  10~(2) (1995) 100--112.
\newblock \href {http://dx.doi.org/10.1007/BF01743537}
  {\path{doi:10.1007/BF01743537}}.

\bibitem{Nguyen-Xuan2017}
H.~Nguyen-Xuan, A polytree-based adaptive polygonal finite element method for
  topology optimization, International Journal for Numerical Methods in
  Engineering 110~(10) (2017) 972--1000, nme.5448.
\newblock \href {http://dx.doi.org/10.1002/nme.5448}
  {\path{doi:10.1002/nme.5448}}.

\bibitem{Wang14ATO}
Y.~Wang, Z.~Kang, Q.~He, Adaptive topology optimization with independent error
  control for separated displacement and density fields, Computers \&
  Structures 135 (2014) 50 -- 61.
\newblock \href {http://dx.doi.org/10.1016/j.compstruc.2014.01.008}
  {\path{doi:10.1016/j.compstruc.2014.01.008}}.

\bibitem{Lambe2018}
A.~B. Lambe, A.~Czekanski, Topology optimization using a continuous density
  field and adaptive mesh refinement, International Journal for Numerical
  Methods in Engineering 113~(3) (2018) 357--373, nme.5617.
\newblock \href {http://dx.doi.org/10.1002/nme.5617}
  {\path{doi:10.1002/nme.5617}}.

\bibitem{Brackett2011}
D.~Brackett, I.~Ashcroft, R.~Hague, Topology optimization for additive
  manufacturing, in: Proceedings of the solid freeform fabrication symposium,
  Austin, TX, Vol.~1, S, 2011, pp. 348--362.

\bibitem{Langelaar16}
M.~Langelaar, An additive manufacturing filter for topology optimization of
  print-ready designs, Struct. Multidiscip. Optim. (2016) 1--13\href
  {http://dx.doi.org/10.1007/s00158-016-1522-2}
  {\path{doi:10.1007/s00158-016-1522-2}}.

\bibitem{Gaynor16}
A.~T. Gaynor, J.~K. Guest, Topology optimization considering overhang
  constraints: Eliminating sacrificial support material in additive
  manufacturing through design, Struct. Multidiscip. Optim. 54~(5) (2016)
  1157--1172.
\newblock \href {http://dx.doi.org/10.1007/s00158-016-1551-x}
  {\path{doi:10.1007/s00158-016-1551-x}}.

\bibitem{Mirzendehdel2016CAD}
A.~M. Mirzendehdel, K.~Suresh, Support structure constrained topology
  optimization for additive manufacturing, Comput. Aided Des. 81~(C) (2016)
  1--13.
\newblock \href {http://dx.doi.org/10.1016/j.cad.2016.08.006}
  {\path{doi:10.1016/j.cad.2016.08.006}}.

\bibitem{Qian2017}
X.~Qian, Undercut and overhang angle control in topology optimization: A
  density gradient based integral approach, International Journal for Numerical
  Methods in Engineering 111~(3) (2017) 247--272, nme.5461.
\newblock \href {http://dx.doi.org/10.1002/nme.5461}
  {\path{doi:10.1002/nme.5461}}.

\bibitem{Lazarov16}
B.~S. Lazarov, F.~Wang, O.~Sigmund, Length scale and manufacturability in
  density-based topology optimization, Archive of Applied Mechanics 86~(1)
  (2016) 189--218.
\newblock \href {http://dx.doi.org/10.1007/s00419-015-1106-4}
  {\path{doi:10.1007/s00419-015-1106-4}}.

\bibitem{Clausen16}
A.~Clausen, N.~Aage, O.~Sigmund, Exploiting additive manufacturing infill in
  topology optimization for improved buckling load, Engineering 2~(2) (2016)
  250 -- 257.
\newblock \href {http://dx.doi.org/10.1016/J.ENG.2016.02.006}
  {\path{doi:10.1016/J.ENG.2016.02.006}}.

\bibitem{Wu17-infill}
J.~Wu, N.~Aage, R.~Westermann, O.~Sigmund, Infill optimization for additive
  manufacturing -- approaching bone-like porous structures, IEEE Transactions
  on Visualization and Computer Graphics 24~(2) (2018) 1127--1140.
\newblock \href {http://dx.doi.org/10.1109/TVCG.2017.2655523}
  {\path{doi:10.1109/TVCG.2017.2655523}}.

\bibitem{Livesu2017}
M.~Livesu, S.~Ellero, J.~Martínez, S.~Lefebvre, M.~Attene, From 3d models to
  3d prints: an overview of the processing pipeline, Computer Graphics Forum
  36~(2) (2017) 537--564.
\newblock \href {http://dx.doi.org/10.1111/cgf.13147}
  {\path{doi:10.1111/cgf.13147}}.

\bibitem{Wang13}
W.~Wang, T.~Y. Wang, Z.~Yang, L.~Liu, X.~Tong, W.~Tong, J.~Deng, F.~Chen,
  X.~Liu, Cost-effective printing of 3d objects with skin-frame structures, ACM
  Trans. Graph. 32~(6) (2013) 177:1--177:10.
\newblock \href {http://dx.doi.org/10.1145/2508363.2508382}
  {\path{doi:10.1145/2508363.2508382}}.

\bibitem{Zhang2015}
X.~Zhang, Y.~Xia, J.~Wang, Z.~Yang, C.~Tu, W.~Wang, Medial axis tree—an
  internal supporting structure for 3d printing, Computer Aided Geometric
  Design 35-36 (2015) 149 -- 162, geometric Modeling and Processing 2015.
\newblock \href {http://dx.doi.org/10.1016/j.cagd.2015.03.012}
  {\path{doi:10.1016/j.cagd.2015.03.012}}.

\bibitem{Jiang2017}
C.~Jiang, C.~Tang, H.-P. Seidel, P.~Wonka, Design and volume optimization of
  space structures, ACM Trans. Graph. 36~(4) (2017) 159:1--159:14.
\newblock \href {http://dx.doi.org/10.1145/3072959.3073619}
  {\path{doi:10.1145/3072959.3073619}}.

\bibitem{Lu14}
L.~Lu, A.~Sharf, H.~Zhao, Y.~Wei, Q.~Fan, X.~Chen, Y.~Savoye, C.~Tu,
  D.~Cohen-Or, B.~Chen, Build-to-last: Strength to weight {3D} printed objects,
  ACM Trans. Graph. 33~(4) (2014) 97:1--97:10.
\newblock \href {http://dx.doi.org/10.1145/2601097.2601168}
  {\path{doi:10.1145/2601097.2601168}}.

\bibitem{Musialski2016}
P.~Musialski, C.~Hafner, F.~Rist, M.~Birsak, M.~Wimmer, L.~Kobbelt, Non-linear
  shape optimization using local subspace projections, ACM Trans. Graph. 35~(4)
  (2016) 87:1--87:13.
\newblock \href {http://dx.doi.org/10.1145/2897824.2925886}
  {\path{doi:10.1145/2897824.2925886}}.

\bibitem{Zhao2017}
H.~Zhao, W.~Xu, K.~Zhou, Y.~Yang, X.~Jin, H.~Wu, Stress-constrained thickness
  optimization for shell object fabrication, Computer Graphics Forum 36~(6)
  (2017) 368--380.
\newblock \href {http://dx.doi.org/10.1111/cgf.12986}
  {\path{doi:10.1111/cgf.12986}}.

\bibitem{Li2017CGF}
W.~Li, A.~Zheng, L.~You, X.~Yang, J.~Zhang, L.~Liu, Rib-reinforced shell
  structure, Computer Graphics Forum 36~(7) (2017) 15--27.
\newblock \href {http://dx.doi.org/10.1111/cgf.13268}
  {\path{doi:10.1111/cgf.13268}}.

\bibitem{Martinez16}
J.~Mart\'{\i}nez, J.~Dumas, S.~Lefebvre, Procedural voronoi foams for additive
  manufacturing, ACM Trans. Graph. 35~(4) (2016) 44:1--44:12.
\newblock \href {http://dx.doi.org/10.1145/2897824.2925922}
  {\path{doi:10.1145/2897824.2925922}}.

\bibitem{Martinez2017}
J.~Mart\'{\i}nez, H.~Song, J.~Dumas, S.~Lefebvre, Orthotropic k-nearest foams
  for additive manufacturing, ACM Trans. Graph. 36~(4) (2017) 121:1--121:12.
\newblock \href {http://dx.doi.org/10.1145/3072959.3073638}
  {\path{doi:10.1145/3072959.3073638}}.

\bibitem{Schumacher15}
C.~Schumacher, B.~Bickel, J.~Rys, S.~Marschner, C.~Daraio, M.~Gross,
  Microstructures to control elasticity in 3d printing, ACM Trans. Graph.
  34~(4) (2015) 136:1--136:13.
\newblock \href {http://dx.doi.org/10.1145/2766926}
  {\path{doi:10.1145/2766926}}.

\bibitem{Panetta15}
J.~Panetta, Q.~Zhou, L.~Malomo, N.~Pietroni, P.~Cignoni, D.~Zorin, Elastic
  textures for additive fabrication, ACM Trans. Graph. 34~(4) (2015)
  135:1--135:12.
\newblock \href {http://dx.doi.org/10.1145/2766937}
  {\path{doi:10.1145/2766937}}.

\bibitem{Chougrani2017}
L.~Chougrani, J.-P. Pernot, P.~Véron, S.~Abed, Lattice structure lightweight
  triangulation for additive manufacturing, Computer-Aided Design 90 (2017) 95
  -- 104, sI:SPM2017.
\newblock \href {http://dx.doi.org/10.1016/j.cad.2017.05.016}
  {\path{doi:10.1016/j.cad.2017.05.016}}.

\bibitem{Zhu2017}
B.~Zhu, M.~Skouras, D.~Chen, W.~Matusik, Two-scale topology optimization with
  microstructures, ACM Trans. Graph. 36~(5) (2017) 164:1--164:16.
\newblock \href {http://dx.doi.org/10.1145/3095815}
  {\path{doi:10.1145/3095815}}.

\bibitem{Wu2017CMAME}
J.~Wu, A.~Clausen, O.~Sigmund, Minimum compliance topology optimization of
  shell-infill composites for additive manufacturing, Computer Methods in
  Applied Mechanics and Engineering 326 (2017) 358--375.
\newblock \href {http://dx.doi.org/10.1016/j.cma.2017.08.018}
  {\path{doi:10.1016/j.cma.2017.08.018}}.

\bibitem{Liu2017}
X.~Liu, V.~Shapiro, Sample-based synthesis of two-scale structures with
  anisotropy, Computer-Aided Design 90 (2017) 199 -- 209, sI:SPM2017.
\newblock \href {http://dx.doi.org/10.1016/j.cad.2017.05.013}
  {\path{doi:10.1016/j.cad.2017.05.013}}.

\bibitem{Gibson10}
I.~Gibson, D.~W. Rosen, B.~Stucker, Design for Additive Manufacturing, Springer
  US, Boston, MA, 2010, pp. 299--332.
\newblock \href {http://dx.doi.org/10.1007/978-1-4419-1120-9_11}
  {\path{doi:10.1007/978-1-4419-1120-9_11}}.

\bibitem{Pottmann2007}
H.~Pottmann, Architectural geometry, Vol.~10, Bentley Institute Press, 2007.

\bibitem{Norato2015}
J.~Norato, B.~Bell, D.~Tortorelli, A geometry projection method for
  continuum-based topology optimization with discrete elements, Computer
  Methods in Applied Mechanics and Engineering 293 (2015) 306 -- 327.
\newblock \href {http://dx.doi.org/10.1016/j.cma.2015.05.005}
  {\path{doi:10.1016/j.cma.2015.05.005}}.

\bibitem{Sigmund01}
O.~Sigmund, A 99 line topology optimization code written in matlab, Struct.
  Multidiscip. Optim. 21~(2) (2001) 120--127.
\newblock \href {http://dx.doi.org/10.1007/s001580050176}
  {\path{doi:10.1007/s001580050176}}.

\bibitem{Wu16}
J.~Wu, C.~Dick, R.~Westermann, A system for high-resolution topology
  optimization, IEEE Transactions on Visualization and Computer Graphics 22~(3)
  (2016) 1195--1208.
\newblock \href {http://dx.doi.org/10.1109/TVCG.2015.2502588}
  {\path{doi:10.1109/TVCG.2015.2502588}}.

\bibitem{Andreassen10}
E.~Andreassen, A.~Clausen, M.~Schevenels, B.~S. Lazarov, O.~Sigmund, Efficient
  topology optimization in matlab using 88 lines of code, Structural and
  Multidisciplinary Optimization 43~(1) (2010) 1--16.
\newblock \href {http://dx.doi.org/10.1007/s00158-010-0594-7}
  {\path{doi:10.1007/s00158-010-0594-7}}.

\bibitem{Guest04}
J.~K. Guest, J.~Pr{\'e}vost, T.~Belytschko, Achieving minimum length scale in
  topology optimization using nodal design variables and projection functions,
  International Journal for Numerical Methods in Engineering 61~(2) (2004)
  238--254.
\newblock \href {http://dx.doi.org/10.1002/nme.1064}
  {\path{doi:10.1002/nme.1064}}.

\bibitem{Wang10}
F.~Wang, B.~S. Lazarov, O.~Sigmund, On projection methods, convergence and
  robust formulations in topology optimization, Structural and
  Multidisciplinary Optimization 43~(6) (2010) 767--784.
\newblock \href {http://dx.doi.org/10.1007/s00158-010-0602-y}
  {\path{doi:10.1007/s00158-010-0602-y}}.

\bibitem{Svanberg87}
K.~Svanberg, The method of moving asymptotes---a new method for structural
  optimization, International Journal for Numerical Methods in Engineering
  24~(2) (1987) 359--373.
\newblock \href {http://dx.doi.org/10.1002/nme.1620240207}
  {\path{doi:10.1002/nme.1620240207}}.

\bibitem{Sigmund07}
O.~Sigmund, Morphology-based black and white filters for topology optimization,
  Structural and Multidisciplinary Optimization 33~(4) (2007) 401--424.
\newblock \href {http://dx.doi.org/10.1007/s00158-006-0087-x}
  {\path{doi:10.1007/s00158-006-0087-x}}.

\end{thebibliography}

\section*{Appendix: Sensitivity Analysis}
\label{sec:Sensitivity}
The sensitivity analysis for gradient-based optimization is presented in the following. By using the chain rule, the sensitivity of the objective function $c$ with respect to the design variable $\bm{x}^k$ is given by:
\begin{equation}
\frac{\partial c}{\partial \bm{x}^k} = {\textstyle{\sum^{\bar{k}}\limits_{l=1}}} \left( \frac{\partial c}{\partial \bm{\rho}} \frac{\partial \bm{\rho}}{\partial \tilde{\bm{x}}^l} \frac{\partial \tilde{\bm{x}}^l}{\partial \bm{x}^k} \right).
\label{eq:dcdx}
\end{equation}

The first term ${\partial c}/{\partial \bm{\rho}}$ is derived from adjoint analysis. The individual entry is
\begin{equation}
\frac{\partial c}{\partial \rho_e} = -p (\rho_e)^{p-1} (E_0 - E_{\text{min}})\bm{u}_e^{\mathrm{T}} \bm{k}_0 \bm{u}_e.
\end{equation}

The second term can be evaluated from Eq.~\eqref{eq:MappingFilter}, 
\begin{equation}
\frac{\partial \bm{\rho}}{\partial \tilde{\bm{x}}^l} = \bm{T}^l.
\end{equation}

The third term comes from the smoothed refinement filter, Eq.~\eqref{eq:filter}. The non-zero entries (i.e., if $l \geq k$) are calculated by 
\begin{equation}
\frac{\partial \tilde{{x}}^l_{i,j}}{\partial {x}^k_{i^{k-l},j^{k-l}}} = \frac{1}{l} \left( \frac{1}{l} \, {\textstyle\sum^{l-1}\limits_{m=0}} \left( x^{l-m}_{i^{-m},j^{-m}} \right)^{p_n} \right)^{\frac{1}{p_n}-1} \left(x^k_{i^{k-l},j^{k-l}}\right)^{p_n - 1}.
\end{equation}

The sensitivity of the volume $V$ with respect to the design variable $\bm{x}^k$ is derived in a similar way:
\begin{equation}
\frac{\partial V}{\partial \bm{x}^k} = {\textstyle{\sum^{\bar{k}}\limits_{l=k}}} \left( \frac{\partial V}{\partial \bm{\rho}} \frac{\partial \bm{\rho}}{\partial \tilde{\bm{x}}^l} \frac{\partial \tilde{\bm{x}}^l}{\partial \bm{x}^k} \right).
\end{equation}
Individual entries in the first term are given by
\begin{equation}
\frac{\partial V}{\partial \rho_e} = v_e,
\end{equation}
where $v_e$ is the unit volume for each element.


\end{document}